\documentclass[psamsfonts]{amsart}

\usepackage{amssymb,amsfonts}
\usepackage{amsmath}
\usepackage[all,arc]{xy}
\usepackage{enumerate}
\usepackage{mathrsfs}
\usepackage{mathtools}
\usepackage{cite}
\usepackage[unicode]{hyperref}
\usepackage{bookmark}

\newtheorem{thm}{Theorem}[section]
\newtheorem{theorem}[thm]{Theorem}
\newtheorem{cor}[thm]{Corollary}
\newtheorem{prop}[thm]{Proposition}
\newtheorem{lemma}[thm]{Lemma}

\newtheorem{quest}[thm]{Question}
\newtheorem{claim}[thm]{Claim}
\newtheorem{fact}[thm]{Fact}

\theoremstyle{definition}
\newtheorem{defn}[thm]{Definition}

\newtheorem{rem}[thm]{Remark}

\newcommand{\defeq}{\mathrel{\mathop:}=}

\newcommand{\Q}{\mathbb{Q}}
\newcommand{\PP}{\mathbb{P}}

\newcommand{\LL}{\mathbb{L}}
\newcommand{\CC}{\mathbb{C}}

\newcommand{\UU}{\mathbb{U}}
\newcommand{\bbS}{\mathbb{S}}
\newcommand{\BB}{\mathbb{B}}
\newcommand{\HH}{\mathbb{H}}
\newcommand{\FF}{\mathbb{F}}
\newcommand{\A}{\mathbb{A}}
\newcommand{\SSS}{\mathbb{S}}
\newcommand{\bbK}{\mathbb{K}}

\newcommand{\la}{\lambda}
\newcommand{\ka}{\kappa}
\newcommand{\w}{\omega}
\newcommand{\sm}{\setminus}
\newcommand{\mc}{\mathcal}

\DeclareMathOperator{\cf}{cf}

\DeclareMathOperator{\crit}{crit}
\DeclareMathOperator{\ITP}{ITP}
\DeclareMathOperator{\Add}{Add}
\DeclareMathOperator{\dom}{dom}

\DeclareMathOperator{\Coll}{Col}
\DeclareMathOperator{\Index}{Index}
\DeclareMathOperator{\Ord}{Ord}

\DeclareTextCommand{\textaleph}{L8U}{ℵ}

\newcommand{\uhr}{\upharpoonright}
\newcommand{\forces}{\Vdash}
\newcommand{\cat}{^\smallfrown}

\makeatletter
\let\c@equation\c@thm
\makeatother
\numberwithin{equation}{section}

\title{The Strong and Super Tree Properties at Successors of Singular Cardinals}
\author{William Adkisson}
\begin{document}

\begin{abstract}
	The strong tree property and ITP (also called the super tree property) are generalizations of the tree property that characterize strong compactness and supercompactness up to inaccessibility. That is, an inaccessible cardinal $\kappa$ is strongly compact if and only if the strong tree property holds at $\kappa$, and supercompact if and only if ITP holds at $\kappa$.
	We present several results motivated by the problem of obtaining the strong tree property and ITP at many successive cardinals simultaneously; these results focus on the successors of singular cardinals.
	We describe a general class of forcings that will obtain the strong tree property and ITP at the successor of a singular cardinal of any cofinality. 
	Generalizing a result of Neeman about the tree property, we show that it is consistent for ITP to hold at $\aleph_n$ for all $2 \leq n < \omega$ simultaneously with the strong tree property at $\aleph_{\omega+1}$; we also show that it is consistent for ITP to hold at $\aleph_n$ for all $3 < n < \omega$ and at $\aleph_{\omega+1}$ simultaneously. Finally, turning our attention to singular cardinals of uncountable cofinality, we show that it is consistent for the strong and super tree properties to hold at successors of singulars of multiple cofinalities simultaneously.
\end{abstract}

	\maketitle

\section{Introduction}
Obtaining compactness properties at many successive cardinals simultaneously is a long-running project of set theory. The tree property and its generalizations have been of particular interest. The starting point for this program was a result of Mitchell and Silver \cite{mitchell:tp}, who showed that it is consistent from a weakly compact cardinal for the tree property to hold at $\aleph_2$; their techniques have formed the blueprint for obtaining the tree property at successors of regular cardinals. Abraham \cite{abraham:treepropN2N3} showed that it is possible for the tree property to hold at $\aleph_2$ and $\aleph_3$ simultaneously, starting from a supercompact cardinal and a weakly compact cardinal. Cummings and Foreman \cite{CF_TreeProp}, starting with $\w$ many supercompacts, produced a model in which the tree property holds at $\aleph_n$ for $2 \leq n < \w$.

Different techniques are required to obtain the tree property at the successor of a singular cardinal.
Magidor and Shelah \cite{MS_TPSuccSing} showed that the tree property could hold at $\aleph_{\w+1}$ using very strong large cardinal assumptions; these assumptions were reduced to $\w$ supercompacts by Sinapova \cite{sinapova:tpatNw+1}. Similar results can be obtained at the successors of singular cardinals of uncountable cofinality. Building off techniques of Neeman \cite{neeman:tpnotschcountablecof}, Sinapova \cite{sinapova:tpnotschunccof} showed that the tree property could consistently hold at the successor of a singular cardinal of any cofinality (along with the failure of SCH). Golshani \cite{golshani:tpunc} further refined this result, showing that the tree property could hold at $\aleph_{\gamma+1}$ for any limit ordinal $\gamma$.
Neeman \cite{NeemanTPNw+1} was able to combine the techniques for obtaining the tree property at successors of regular cardinals and successors of singular cardinals, producing a model in which the tree property holds at $\aleph_n$ for $2\leq n < \w$ and at $\aleph_{\w+1}$ simultaneously.

We now turn our attention to generalizations of the tree property. Two generalizations of interest, due to Jech and Magidor respectively, are the strong tree property and the super tree property (also known as ITP). These properties characterize strong compactness and supercompactness up to inaccessibility in much the same way that the tree property characterizes weak compactness. In particular, if $\ka$ is inaccessible, then:
\begin{itemize}
	\item The tree property holds at $\ka$ if and only if $\ka$ is weakly compact
	\item (Jech \cite{jech:combinatorialprobs}) The strong tree property holds at $\ka$ if and only if $\ka$ is strongly compact
	\item (Magidor \cite{magidor:combinatorialcharsc}) The super tree property (ITP) holds at $\ka$ if and only if $\ka$ is supercompact.
\end{itemize}
Although these properties were originally examined in the 1970s, recent work of Wei\ss, who studied these properties in his thesis \cite{weiss:thesis}, has sparked a number of new results. Fontanella \cite{fontanella:ITPN2N3} generalized Abraham's result to obtain ITP at $\aleph_2$ and $\aleph_3$ simultaneously. Shortly thereafter, Fontanella \cite{fontanella:CF} and Unger \cite{UngerCF} independently showed that in the model of Cummings and Foreman, ITP holds at each $\aleph_n$. Fontanella \cite{fontanella:stpsuccsing} also showed that assuming $\w$-many supercompact cardinals, the strong tree property could consistently hold at $\aleph_{\w+1}$. 
Hachtman and Sinapova \cite{HachtmanITPNw+1} generalized this result further, showing that from the same assumptions it is consistent to have ITP at $\aleph_{\w+1}$.

In this paper, we continue the analysis of the strong and super tree properties at successors of singular cardinals and at small successor cardinals simultaneously. In Section \ref{s:TPuptoNw+1}, we show that in Neeman's construction in \cite{NeemanTPNw+1}, ITP holds at each $\aleph_n$, and the strong tree property holds at $\aleph_{\w+1}$. This answers a question of Fontanella, who asked if it was consistent to have the strong tree property up to $\aleph_{\w+1}$. The key step in this argument is Theorem \ref{thm:TP2cardfullNA}, 
which describes a general class of forcings in which the strong tree property holds at the successor of a singular. This theorem applies to singular cardinals of any cofinality.

Next, we consider ITP at successors of singular cardinals. In Section \ref{s:ITPsuccsing}, we describe a similar class of forcings in which ITP holds at the successor of a singular. Once again, there are no restrictions on the cofinality of this singular cardinal. We use this to show that from $\tau$-many supercompacts it is consistent to have ITP at $\aleph_{\tau+1}$ for any regular $\tau$ with $\tau < \aleph_\tau$. In Section \ref{s:ITPalmostuptoNw+1}, we apply this technique to show that in a modified version of Neeman's model, ITP holds at $\aleph_{\w+1}$ and at $\aleph_n$ for all $n > 3$. This is a partial answer to a question of Hachtman and Sinapova, who asked if ITP could be obtained up to $\aleph_{\w+1}$.

Finally, in Section \ref{s:manycofs} we show that these techniques can also be applied to obtain the strong and/or super tree properties at successors of small singular cardinals of multiple cofinalities simultaneously. In particular, given a finite increasing sequence $\alpha_0, \dots, \alpha_n$ of ordinals where for each $i \leq n$ $\alpha_i < \omega_{\alpha_i}$ and $\omega_{\alpha_i}$ is a regular cardinal, we describe a model in which we obtain these properties at $\aleph_{\omega_{\alpha_i}+1}$ for all $i \leq n$.
If $\alpha_{i+1} > \alpha_i+1$, we obtain ITP at $\aleph_{\omega_{\alpha_i}+1}$; otherwise, we can obtain the strong tree property.
Our large cardinal hypotheses here are somewhat stronger than in other sections; we require $\ka_0$-many supercompacts, where $\ka_0$ is itself supercompact.

On notation: we use the convention that a poset is $\ka$-closed if every descending sequence of size less than $\ka$ has a lower bound. If $\PP$ is a forcing notion over a model $V$, we will use $V[\PP]$ to denote the extension of $V$ by some generic for $\PP$.

\section{Preliminaries and Branch Lemmas}\label{section:lemmas}

We begin by defining the the strong tree property and ITP.

\begin{defn}
	Let $\mu$ be a regular cardinal and let $\la \geq \mu$. A set $d = \langle d_z \mid z\in \mc{P}_\mu(\la)\rangle$ is a \emph{$\mc{P}_\mu(\la)$-list} if for all $z \in \mc{P}_\mu(\la)$, $d_z\subseteq z$. 
	We define the $z$-th level of a list, denoted $L_z$, by
	\[L_z = \{d_y \cap z \mid z\subseteq y, y \in \mc{P}_\mu(\la)\}.\]
	A $\mc{P}_\mu(\la)$-list is \emph{thin} if $|L_z| < \mu$ for all $z \in \mc{P}_\mu(\la)$.
\end{defn}

Note that if $\mu$ is inaccessible, every $\mc{P}_\mu(\la)$ list is thin.

\begin{defn}
	A set $b \subseteq \la$ is a cofinal branch through a $\mc{P}_\mu(\la)$-list $d$ if for all $x \in \mc{P}_\mu(\la)$ there is $z \supseteq x$ such that $b\cap x = d_z\cap x$. In other words, for all $x$, $b \cap x \in L_x$. We say that the strong tree property holds at $\mu$ if for all $\lambda \geq \mu$, every thin $\mc{P}_\mu(\la)$-list has a cofinal branch.
\end{defn}

\begin{defn}
	A set $b \subseteq \la$ is an ineffable branch through a $\mc{P}_\mu(\la)$-list $d$ if $\{z \in \mc{P}_\mu(\la) \mid b \cap z = d_z\}$ is stationary. We say $\ITP(\mu,\la)$ holds if every thin $\mc{P}_\mu(\la)$-list has an ineffable branch. We say ITP holds at $\mu$ if $\ITP(\mu,\la)$ holds for all $\la \geq \mu$.
\end{defn}

Note that every ineffable branch is cofinal.

\begin{fact}\cite[Proposition 3.4]{weiss:combinatorialessence}
	Let $\la' > \la$. $\ITP(\mu, \la')$ implies $\ITP(\mu, \la)$.
\end{fact}

\begin{fact}\cite[Lemma 3.4]{fontanella:stpsuccsing}\label{fact:TPcandecreaselambda}
	Let $\la' > \la$. If every thin $\mc{P}_\mu(\la')$-list has a cofinal branch, then so does every thin $\mc{P}_\mu(\la)$ list.
\end{fact}

ITP, like the tree property, is usually obtained at the successor of regular cardinals by means of a lifted embedding. This will produce a branch in the generic extension containing the embedding, as described in the following lemma.

\begin{lemma}\label{lem:embeddingtobranch}
	Let $W$ be a model of set theory, and let $d$ be a thin $\mc{P}_\ka(\la)$ list. Suppose that in some extension $W[G]$ there is a generic embedding $j:W \to M$ with critical point $\ka$ such that $j(\ka) > \la$ and $M^\la \subseteq M$. Then in $W[G]$, $d$ has a cofinal branch $b$. Moreover, if $b \in W$, then $b$ is ineffable.
\end{lemma}
\begin{proof}
	Work in $W[G]$, and consider $j(d)$. Since $j(\ka) > \la$, we have that $j''\la \in \mathcal{P}_{j(\ka)}(j(\la))$.
	Let $b = \{\alpha < \la \mid j(\alpha) \in j(d)_{j''\la}\}$.
	
	We first claim that $b$ is a cofinal branch. Let $z \in \mathcal{P}_\ka(\la)^W$. We want to show that $b \cap z \in L_z$. Note that $j(d)_{j''\la}\cap j(z) \in j(L_z)$.
	Since $d$ is a thin list, $|L_z|< \ka$, so $j(L_z) = j''(L_z) = \{j(d)_{j(y)} \cap j(z) \mid z\subseteq y, y \in \mc{P}_\ka(\la)^W\}$. Since $j(d)_{j''\la}\cap j(z) \in j(L_z)$, there must be some $y$ such that $j(d)_{j''\la}\cap j(z) = j(d)_{j(y)}\cap j(z)$. It follows that $b\cap z = d_y\cap z$.
	
	Now, suppose that $b \in W$. Let $U$ be the normal measure on $\mc{P}_\ka(\la)^W$ corresponding to $j$. Note that in the ultrapower of $W$ by $U$, $j''\la$ is represented by $[x\to x]$, and so $j(d)_{j''\la}$ is represented by $[x \mapsto d]_{[x\mapsto x]} = [x \mapsto d_x]$. On the other hand, since $b \in W$, we can take $j(b)$; indeed, $j(d)_{j''\la} = j(b) \cap j''\la$. Since $j(b)$ is represented by $[x \mapsto b]$, $j(d)_{j''\la} = [x \mapsto b] \cap [x \mapsto x]$. We conclude that $d_x = b \cap x$ for $U$-many $x$. Since $U$ is a normal measure on $W$, all measure one sets must be stationary in $W$, so $b$ must be an ineffable branch.
\end{proof}

The next step is usually to show that the forcing could not have added a branch, so the branch must have already been present in $W$, and will thus be ineffable.
The crucial ingredients in this process are branch preservation lemmas. To generalize the standard branch lemmas from trees to thin lists, we use the (thin) approximation property.

\begin{defn}\label{def:ka-approx}
	Let $\ka$ be regular. A forcing $\PP$ has the \emph{$\ka$-approximation property} in a model $V$ if for every ordinal $\la$ and every $\PP$-name $\dot{b}$ for a subset of $\la$, $\forces_\PP \dot{b} \cap z \in V$ for all $z \in (\mc{P}_\ka(\la))^V$ implies that $\forces_\PP \dot{b} \in V.$

	A poset $\PP$ has the \emph{thin $\ka$-approximation property} in $V$ if for every ordinal $\la$ and every $\PP$-name $\dot{b}$ for a subset of $\la$, $\forces_\PP \dot{b} \cap z \in V$ and $|\{x \in V \mid \exists p\in \PP \ p \forces_\PP x = \dot{b}\cap z\}|<\ka$ for every $z \in (\mc{P}_\ka(\la))^V$ implies that $\forces_\PP \dot{b} \in V$.
\end{defn}

Note that the thin $\ka$-approximation property is weaker than the $\ka$-approximation property. We say that names meeting the hypotheses in Definition \ref{def:ka-approx} are (thinly) $\ka$-approximated by $\PP$ over $V$. 
These properties are useful because a cofinal branch through a thin $\mc{P}_\ka(\la)$ list is always thinly approximated over $V$.

\begin{lemma}\label{branchapprox}
	Let $d$ be a thin $\mc{P}_\ka(\la)$ list in $V$, and let $\mathbb{P}$ be a notion of forcing over $V$. Suppose $\dot{b}$ is a $\PP$-name for a cofinal branch through this list. Then $\dot{b}$ is thinly $\ka$-approximated by $\PP$ over $V$.
\end{lemma}
\begin{proof}
	Let $\dot{b}$ be a name for a cofinal branch, and let $x \in (\mc{P}_\ka(\la))^V$. Since $\dot{b}$ is forced to be cofinal, $\dot{b}\cap x$ is forced to be in $L_x$ for all $x$. Since $L_x$ is in $V$, $\dot{b}\cap x$ is likewise forced to be in $V$. Since the list is thin, $|L_x| < \ka$, and so there are fewer than $\ka$ possibilities for $\dot{b}\cap x$.
\end{proof}

\begin{lemma}\label{lem:cctoapprox}\cite[Lemma 7.10]{UngerCF}
	Suppose that $\dot{b}$ is a $\PP\ast\dot{\Q}$-name for a subset of some ordinal $\mu$, which is (thinly) $\ka$-approximated by $\PP\ast \dot{\Q}$ over $V$. If $\PP$ has the $\ka$-cc, then in $V[\PP]$, $\dot{b}$ is still (thinly) $\ka$-approximated by $\Q$.
\end{lemma}
\begin{proof}
	Let $G$ be generic for $\PP$ over $V$. Since $\PP$ is $\ka$-cc, $\mc{P}_\ka(\mu)^V$ is cofinal in $\mc{P}_\ka(\mu)^{V[G]}$. In $V[G]$, if $x \in \mc{P}_\ka(\mu)$ is a subset of some $y \in \mc{P}_\ka(\mu)^{V}$, then $\dot{b}\cap x$ is determined solely by $\dot{b}\cap y$ and $x$. Since $\dot{b}$ is $\ka$-approximated by $\PP\ast\dot{\Q}$ over $V$, $\dot{b}\cap y$ must be in $V$. It follows that $\Q$ forces that $\dot{b}\cap x \in V[G]$, as desired.
	
	Now suppose that $\dot{b}$ is thinly $\ka$-approximated by $\PP\ast\dot{\Q}$, not just $\ka$-approximated. For every possible value of $\dot{b}\cap x$, we can always choose an extension that is a possible value for $\dot{b}\cap y$. This defines an injection from the possible values (in $V[G]$) of $\dot{b}\cap x$ to the possible values (in $V)$ of $\dot{b}\cap y$. Since $\dot{b}$ is thinly $\ka$-approximated in $V$ by $\PP\ast\dot{\Q}$, there $<\ka$-many possible values for $\dot{b}\cap y$, and thus there will be $<\ka$-many possible values in $V[G]$ for $\dot{b}\cap x$. So $\dot{b}$ is thinly $\ka$-approximated by $\Q$ over $V[G]$.
\end{proof}

\begin{lemma}\label{lem:ccpreservesdistributive}(Part of Easton's Lemma)
	Let $\PP$ be $\ka$-closed in $W$, where $V$ is a $\ka$-cc forcing extension of $W$. Then forcing with $\PP$ over $V$ does not add any sequences of length $<\ka$.
\end{lemma}

\begin{lemma}\label{lem:cc+closedtoapprox}
	Suppose that $\dot{b}$ is a $\PP \ast \dot{\Q} \ast \dot{\SSS}$-name for a subset of some ordinal $\mu$, which is $\ka$-approximated. If $\PP$ is $\ka$-cc, and $\Q$ does not add sequences of length less than $\ka$ to $V[\PP]$, then in $V[\PP \ast \dot{\Q}]$, $\dot{b}$ is still $\ka$-approximated by $\SSS$. 
\end{lemma}
\begin{proof}
	Let $z \in (\mc{P}_\ka(\mu))^{V[\PP\ast \dot{\Q}]}$. Since $\Q$ does not add sequences of length $< \ka$, $z \in (\mc{P}_\ka(\mu))^{V[\PP]}$. By Lemma \ref{lem:cctoapprox}, $\forces_{\Q\ast\dot{\SSS}} z \cap \dot{b} \in V[\PP]$. We conclude that, working in $V[\PP \ast \dot{\Q}]$, $\forces_{\SSS} z \cap \dot{b} \in V[\PP] \subseteq V[\PP\ast \dot{\Q}]$.
\end{proof}

To show that forcings have the (thin) approximation property, we will use the following lemmas, which generalize traditional branch preservation lemmas for trees.

\begin{lemma}\label{lem:closedbranch}\cite[Lemma 7.9]{UngerCF}
	Let $\tau < \ka$ with $2^\tau \geq \ka$. If $\PP$ is $\tau^+$-closed, then $\PP$ has the thin $\ka$-approximation property.
\end{lemma}

\begin{lemma}\label{lem:ccbranch}\cite[Lemma 7.8]{UngerCF}
	Let $\ka$ be regular. Suppose $\PP$ is a poset such that $\PP\times \PP$ is $\ka$-cc. Then $\PP$ has the $\ka$-approximation property.
\end{lemma}

The following definition and lemma provide a tool for verifying when the square of a poset has the $\ka$-cc in outer models.

\begin{defn}
	Let $K\subseteq V$ be a model of a sufficiently large fragment of ZFC. We say that $K$ has the $<\delta$-covering property with respect to $V$ if for every $A \subseteq K$ in $V$ with $|A|<\delta$, there is some $B \in K$ so that $K\models |B|<\delta$ and $A \subseteq B$.
\end{defn}

\begin{lemma}\label{lem:covtocc}\cite[Claim 2.2]{NeemanTPNw+1}
	Let $\delta < \ka$ be regular cardinals, and suppose $K$ is a model of some large enough fragment of ZFC that has the $<\ka$-covering property with respect to $V$. Suppose also that for all $\gamma < \ka$, $K \models \gamma^{<\delta} < \ka$. Let $\PP$ be a forcing notion in $K$ whose conditions are functions of size $<\delta$ in $K$. Then any family of size $\ka$ in $V$ of conditions in $\PP$ can be refined to a family of the same size whose domains form a $\Delta$-system.
\end{lemma}

We will also use the following lemma that generalizes \cite[Claim 2.4]{NeemanTPNw+1}.
\begin{lemma}\label{lem:approxcc+closed}
	Let $\PP$ be a $\ka$-cc forcing and $\Q$ be a $\ka$-closed forcing. Suppose there exists $\tau < \ka$ such that $2^\tau > \ka$. Then $\Q$ has the thin $\ka^+$-approximation property over the generic extension by $\PP$.
\end{lemma}	
\begin{proof}
	Let $\la>\ka$ be some ordinal, and let $\dot{d}$ be a $\PP\times \Q$-name for a subset of $\la$ that is thinly $\ka^+$-approximated by $\Q$ over the generic extension of $V$ by $\PP$. Suppose that $1 \forces_{\PP\times \Q} \dot{d}\notin V[\dot{G}_\PP]$, where $\dot{G}_{\PP}$ is the canonical name for a generic of $\PP$.
	Let $\delta \leq \tau$ be the least such that $2^\delta > \ka$. For each $\sigma \in 2^{<\delta}$ we will build a triple $(A_\sigma, q_\sigma, x_\sigma)$ where $A_\sigma$ is a maximal antichain in $\PP$, $q_\sigma \in \Q$, and $x_\sigma \in (\mc{P}_{\ka^+}(\la))^V$, such that $q_{\sigma'} \leq q_\sigma$ whenever $\sigma'$ extends $\sigma$, and for all $\sigma \in 2^{<\delta}$ and all $p \in A_\sigma$, $(p,q_{\sigma\cat 0})$ and $(p,q_{\sigma\cat 1})$ force contradictory values for $\dot{d}\cap x_\sigma$. Note that it suffices to consider $x_\sigma \in (\mc{P}_{\ka^+}(\la))^V$, since $(\mc{P}_{\ka^+}(\la))^V$ is cofinal in $(\mc{P}_{\ka^+}(\la))^{V[\PP]}$ by the $\ka$-cc of $\PP$.
	
	We will build this inductively, using the following claim.
	
	\begin{claim}\label{cl:splitting1}
		For any $p \in \PP$ and $q_0, q_1\in \Q$, we can find $p' \leq p$, $q_0'\leq q_0$, $q_1'\leq q_1$, and $x \in \mc{P}_\ka(\la)$ such that $p'$ forces that $q'_0$ and $q'_1$ force contradictory information about $\dot{d}\cap x$.
	\end{claim}
	\begin{proof}
	Recall that $1 \forces_{\PP\times\Q} \dot{d}\notin V[\dot{G}_\PP].$  It follows that if $G\times Q_0 \times Q_1$ is generic for $\PP\times \Q \times \Q$, the interpretations $\dot{d}_{G\times Q_0}$ and $\dot{d}_{G\times Q_1}$ must not be equal. Thus there must be a pair of conditions $(p'_0,q_0') \leq (p,q_0)$ and $(p'_1, q'_1) \leq (p,q_1)$ such that $(p'_0,q_0') \in G\times Q_0$ and $(p'_1,q_1') \in G\times Q_1$ and the pair decides incompatible information about $\dot{d}\cap x$ for some $x \in \mc{P}_{\ka^+}(\la)$. Since $p'_0$ and $p'_1$ are both in $G$, they are compatible. We set $p'$ to be a common extension of $p'_0$ and $p'_1$; the conclusion follows.
	\end{proof}

	Let $q \in \Q$. We wish to build $(A,q_0, q_1, x)$ where $A$ is a maximal antichain in $\PP$, $q_0, q_1 \leq q$, and for all $p \in A$, $p$ forces that $q_0$ and $q_1$ force contradictory information about $\dot{d}\cap x$. 
	
	For each $q \in \Q$, we construct $(A, q_0, q_1, x)$ recursively. Set $q_0^0 = q_0^1 = q$ and $A_0 = \emptyset$.
	At successor stages $\alpha+1$ we choose a new $p$ incompatible with all conditions in our antichain-in-construction $A_\alpha$ if such exists, then strengthen it to $p' \leq p$ so that there exist some $q^0_{\alpha+1} \leq q^0_\alpha$, $q^1_{\alpha+1} \leq q^1_\alpha$, and $x_{\alpha+1}\supset x_\alpha$ as desired. We then define $A_{\alpha+1}$ to be $A_\alpha \cup \{p'\}$. At limit stages $\gamma$, we take unions (for $A_\gamma$ and $x_\gamma$) and lower bounds for $q^0_\gamma$ and $q_\gamma^1$. Note that since $\PP$ is $\ka$-cc, this process will end at some $\gamma < \ka$, and so since $\Q$ is $\ka$-closed we will always be able to take these lower bounds. This process gives $(A, q_0, q_1, x)$.
	Applying this construction repeatedly, we can build the binary tree structure as desired.
	
	Now let $x = \bigcup_{\sigma \in 2^{<\delta}} x_\sigma.$ For each $f \in 2^\delta$, by closure we can find $q'_f$ such that $q'_f \leq q_{f\upharpoonright\alpha}$ for all $\alpha < \delta$. Let $G$ be generic for $\PP$. Since $G$ is generic, $G \cap A_\sigma$ is nonempty for all $\sigma$. Working in $V[G]$, we define $q_f \leq q'_f$ to be a refinement deciding $\dot{d}\cap x$. 
	
	We claim that if $f \neq g$ are elements of $2^\delta$, then $q_f$ and $q_g$ force contradictory things about $\dot{d}\cap x$. Let $\alpha$ be the largest such that $f\upharpoonright\alpha = g\upharpoonright\alpha$. Since $G \cap A_{f\uhr \alpha}$ is nonempty, $q_{f\uhr \alpha \cat 0}$ and $q_{f\uhr \alpha \cat 1}$ force contradictory information about $x_{f\uhr \alpha}\subseteq x$. Since $f$ extends one of these and $g$ extends the other, we see that $q_f$ must disagree with $q_g$ on $\dot{d}\cap x$. Thus every $q_f$ for $f\in 2^\delta$ forces a different value for $\dot{d}\cap x$, and so there must be at least $2^\delta > \ka$ potential values for $\dot{d}\cap x$. But since $\dot{d}$ is thinly $\ka^+$-approximated, there are fewer than $\ka^+$ potential values for $\dot{d}$, and thus for $\dot{d}\cap x$. This gives a contradiction.
\end{proof}
At the successor of a singular strong limit cardinal, we have a similar lemma, a generalization of \cite[Lemma 2.1]{MS_TPSuccSing}. The main difference is that we index the construction by $\nu^{<\tau}$ instead of $2^{<\delta}$. To prove this lemma we need a couple of claims, which expand on the splitting behavior described in Claim \ref{cl:splitting1}

\begin{claim}\label{claim:singsplitting1}
	Suppose $\nu$ is a singular strong limit cardinal with cofinality $\tau$, and let $\la \geq \nu^+$. Let $\bbS$ be a notion of forcing. Let $\dot{d}$ be a $\bbS$-name for a subset of $\la$ that is thinly $\nu^+$-approximated by $\bbS$ over $V$, such that $\dot{d}$ is forced not to be in $V$.
	Then for any $s \in \bbS$, there is a sequence $\langle s_\alpha \mid \alpha < \nu\rangle$ of conditions in $\bbS$ and a set $x \in \mc{P}_{\nu^+}(\la)$, such that $s_\alpha \leq s$ for each $\alpha < \nu$ and each $s_\alpha$ decides a different value for $\dot{d}\cap x$.
\end{claim}
\begin{proof}
	Fix $s \in \bbS$. For $x \in \mc{P}_{\nu^+}(\la)$, let $T(x)$ be the set of all $y \subseteq x$ such that some $s' \leq s$ forces $\dot{d}\cap x = y$. Since $\dot{d}$ is thinly $\nu^+$-approximated, $|T(x)| \leq \nu$ for all $x$. We wish to construct $x \in \mc{P}_{\nu^+}(\la)$ such that $|T(x)| = \nu$.
	
	Fix $x_0 \in \mc{P}_{\nu^+}(\la)$. We will build a $\subseteq$-increasing sequence $\langle x_\alpha \mid \alpha < \nu\rangle$ such that for all $\alpha < \nu$ and all $y \in T(x_\alpha)$, there are distinct $y_0$ and $y_1$ in $T(x_{\alpha+1})$ such that $y_0 \cap x_\alpha = y_1 \cap x_\alpha = y$. We do so recursively.
	
	Fix $\alpha < \nu$, and let $y \in T(x_\alpha)$ as witnessed by a condition $s' \leq s$. Since $\dot{d}$ is forced not to be in $V$, we can find a set $x_{\alpha+1}^y \supseteq x_\alpha$, distinct $y'_0, y'_1 \in \mc{P}_{\nu^+}(\la)$, and conditions $s_0, s_1 \leq s'$ such that $s_0 \forces \dot{d}\cap x_{\alpha+1}^y = y'_0$ and $s_1 \forces \dot{d}\cap x_{\alpha+1}^y = y'_1$. (If such objects did not exist, then a single condition could decide $\dot{d}$, contradicting the assumption that $\dot{d} \notin V$.)
	Since $s_0$ and $s_1$ extend $s'$, we must have $y'_0 \cap x_\alpha = y'_1 \cap x_\alpha = y$.
	Let $x_{\alpha+1} = \bigcup_{y \in T(x_\alpha)} x_{\alpha+1}^y$. For each $y \in T(x_\alpha)$, let $s_0$, $s_1$, $y'_0$, and $y'_1$ be as above. Extend $s_0$  to decide the value of $\dot{d}\cap x_{\alpha+1}$ to be some $y_0 \supseteq y'_0$; similarly, extend $s_1$ to force $\dot{d}\cap x_{\alpha+1} = y_1$ for some $y_1 \supseteq y'_1$. Then $y_0$ and $y_1$ are distinct elements of $T(x_{\alpha+1})$, and $y_0\cap x_\alpha = y_1\cap x_\alpha = y$, so $x_{\alpha+1}$ has the desired properties.
	This covers successor stages; at limit stages we take unions.
	
	Now let $x = \bigcup_{\alpha < \nu} x_\alpha$. Let $y \in T(x)$. Then for all $\alpha < \nu$, $y\cap x_\alpha \in T(x_\alpha)$. By construction, we have distinct $y_0$ and $y_1$ in $T(x_\alpha+1)$ (as witnessed by conditions $s_0$ and $s_1$) such that $y_0 \cap x_\alpha = y_1 \cap x_\alpha = y\cap x_\alpha$. Moreover, any extensions of these conditions that decide $\dot{d}\cap x$ must decide incompatible values. Since there are $\nu$-many of these splittings, there must be at least $\nu$ possible values for $\dot{d}\cap x$. We conclude that $|T(x)| = \nu$.
\end{proof}

\begin{claim}\label{cl:singsplitting2}
	Suppose $\nu$ is a singular strong limit cardinal with cofinality $\tau$, and let $\la \geq \nu^+$. Let $\PP$ and $\Q$ be notions of forcing such that $|\PP|\leq \mu$ and $\Q$ is $\mu^+$-closed for some $\mu$ with $\tau < \mu < \nu$. Let $\dot{d}$ be a $\PP\times\Q$-name for a subset of $\la$ that is thinly $\nu^+$-approximated by $\Q$ over the generic extension of $V$ by $\PP$, such that $1_{\PP\times\Q} \forces \dot{d}\notin V[\dot{G}_\PP]$.
	
	Suppose $p_0 \in \PP$, and let $\langle q_\alpha \mid \alpha < \nu\rangle$ be a sequence of conditions in $\Q$. Then there is a sequence $\langle q'_\alpha \mid \alpha < \nu\rangle$ with $q'_\alpha \leq q_\alpha$ for each $\alpha < \nu$, along with  a set $x \in \mc{P}_{\nu^+}(\la)$, such that for all $\alpha < \beta$ there is some $p_{\alpha\beta} \leq p_0$ so that $(p_{\alpha\beta},q'_\alpha)$ and $(p_{\alpha\beta}, q'_\beta$) force incompatible information about $\dot{d}\cap x$.
\end{claim}
\begin{proof}
	Let $G$ be generic for $\PP$ containing $p_0$, and work for the moment in $V[G]$. For all $\alpha < \nu$, by Claim \ref{claim:singsplitting1} there is $x_\alpha \in \mc{P}_{\nu^+}(\la)$ such that there are $\nu$-many possible values for $\dot{d}\cap x_\alpha$ that can be forced by extensions of $q_\alpha$. Returning to $V$, we note that these possible values are also possible values of $\dot{d}\cap x_\alpha$ that can be forced by extensions of $(p_0, q_\alpha)$. Since $\PP$ does not collapse $\nu$ or $\nu^+$, we conclude that (in $V$) there are $\nu$-many possibilities for $\dot{d}\cap x_\alpha$ that can be forced by extensions of $(p_0, q_\alpha)$.
	
	Let $x = \bigcup_{\alpha < \nu} x_\alpha$; note that $x \in \mc{P}_{\nu^+}(\la)$. Since $(p_0,q_\alpha)$ has $\nu$-many extensions deciding pairwise-distinct values for $\dot{d}\cap x_\alpha$, by extending each further we can obtain $\nu$-many extensions of each $(p_0, q_\alpha)$ deciding pairwise-distinct values for $\dot{d}\cap x$. (To do this, simply take the $\nu$-many extensions that disagree on $\dot{d}\cap x_\alpha$, and extend them further to decide $\dot{d}\cap x$. Since $x_\alpha \subseteq x$, these further extensions must also disagree on $\dot{d}\cap x$.)
	
	For all $\alpha < \nu$, we will inductively define a condition $q'_\alpha \leq q_\alpha$ and functions $e_\alpha:\PP\to\PP$ and $f_\alpha:\PP\to \mc{P}(x)$ such that the following hold:
	\begin{itemize}
		\item $e_\alpha(p) \leq p$
		\item $(e_\alpha(p), q'_\alpha)\forces_{\PP\times\Q} \dot{d}\cap x = f_\alpha(p)$
		\item for all $\alpha < \beta < \nu$, $f_\alpha(p) \neq f_\beta(p_0)$ for all $p\in \PP$.
	\end{itemize}
	The sequence $\langle q'_\alpha \mid \alpha < \nu \rangle$ will be the desired object; for each pair $\alpha < \beta$ we will choose $p_{\alpha\beta}$ to be $e_\alpha(e_\beta(p_0))$. We build these objects with a double induction, inducting over both $\alpha$ and the elements of $\PP$.

	Enumerate $\PP$ by $\langle p_\xi \mid \xi < \mu\rangle$, starting with the already determined condition $p_0$. Suppose we have defined $q'_\delta, e_\delta,$ and $f_\delta$ for all $\delta < \alpha$. Inducting on $\xi$, we wish to define $q_\alpha^\xi \leq q_\alpha$, $e_\alpha(p_\xi) \leq p_\xi$, and $f_\alpha(p_\xi)$ such that:
	\begin{itemize}
		\item $\langle q_\alpha^\xi \mid \xi < \mu\rangle$ forms a decreasing sequence below $q_\alpha$,
		\item $(e_\alpha(p_\xi), q_\alpha^\xi)\forces_{\PP\times\Q}\dot{d}\cap x = f_\alpha(p_\xi)$, and
		\item $f_\delta(p) \neq f_\alpha(p_0)$ for all $\delta < \alpha$ and all $p \in \PP$.
	\end{itemize}
	Since $(p_0, q)$ has $\nu$-many extensions deciding $\dot{d}\cap x$ differently, and we have used at most $\max(\alpha, |\PP|) < \nu$ of these extensions in earlier stages, we can pick $e_\alpha(p_0), q_\alpha^0,$ and $f_\alpha(p_0)$ as desired. At successor stages $\xi+1$, we choose an extension $(e_\alpha(p_{\xi+1}), q_\alpha^{\xi+1})$ of $(p_{\xi+1}, q_\alpha^\xi)$ deciding the value of $\dot{d}\cap x$, and set $f_\alpha(p_{\xi+1})$ to be that value. At limit stages $\gamma$, using the $\mu^+$-closure of $\Q$ we select $e_\alpha(p_\gamma) \leq p_\gamma$ and a lower bound $q_\alpha^{\gamma}$ of $\langle q_\alpha^\xi \mid \xi < \gamma \rangle$, such that $(e_\alpha(p_\gamma), q_\alpha^\gamma)$ decides the value of $\dot{d}\cap x$. Set $f_\alpha(p_{\gamma})$ to be this value.
	When this induction on the elements of $\PP$ is finished, we define $q'_\alpha$ to be a lower bound of $\langle q_\alpha^\xi \mid \xi < \mu\rangle$.
	This process constructs functions $e_\alpha$ and $f_\alpha$, as well as the sequence $\langle q'_\alpha\mid \alpha < \nu \rangle$, with the desired properties.
	
	Finally we verify that the sequence $\langle q'_\alpha \mid \alpha < \nu\rangle$ meets the requirements in the statement of the lemma. Fix $\alpha < \beta < \nu$. Let $p_{\alpha\beta} = e_\alpha(e_\beta(p_0))$. By construction, $e_\beta(p_0)\forces \dot{d}\cap x = f_\beta(p_0)$; since $p_{\alpha\beta} \leq e_\beta(p_0)$, $(p_{\alpha\beta}, q'_\beta)$ also forces this. On the other hand, $(p_{\alpha\beta}, q'_\alpha)\forces \dot{d}\cap x = f_\alpha(e_\beta(p_0))$. By construction, these must be different values.
\end{proof}

\begin{lemma}\label{lem:M-S}
Suppose $\nu$ is a singular strong limit cardinal with cofinality $\tau$. Let $\Q$ be a $\mu^+$-closed forcing over a model $V$ for some $\mu < \nu$ with $\tau \leq \mu$, and let $\PP \in V$ be a poset with $|\PP|\leq \mu$. Then $\Q$ has the thin $\nu^+$-approximation property in the generic extension of $V$ by $\PP$.
\end{lemma}
\begin{proof}
	Let $\la$ be a some ordinal, and let $\dot{d}$ be a $\PP\times \Q$-name for a subset of $\la$ that is thinly $\nu^+$-approximated by $\Q$ over the generic extension by $\PP$. Suppose further that the empty condition in $\PP\times \Q$ forces $\dot{d}\notin V[\dot{G}_{\PP}]$. For each $\sigma \in \nu^{<\tau}$ we will build a pair $(q_\sigma, x_\sigma)$ with $q_\sigma \in \Q$ and $x_\sigma \in (\mc{P}_{\nu^+}(\la))^V$, along with a collection of dense sets $\langle D_\sigma^{\alpha\beta}\mid \alpha < \beta < \nu\rangle$, such that $q_{\sigma'} \leq q_\sigma$ whenever $\sigma'$ extends $\sigma$, and for all $\sigma \in \nu^{<\tau}$, all $\alpha < \beta$, and all $p \in D^{\alpha\beta}_{\sigma}$, $(p, q_{\sigma\cat\alpha})$ and $(p, q_{\sigma\cat\beta})$ force contradictory information about $\dot{d}\cap x_{\sigma}$.
	
	Fix $q \in \Q$. We will construct $\langle q_\alpha \mid \alpha < \nu\rangle, x \in \mc{P}_{\nu^+}(\la)$, and a dense sets $D_{\alpha\beta}$, such that each $q_\alpha \leq q$, and for all $\alpha < \beta$ and all $p \in D_{\alpha\beta}$, $p$ forces that $q_\alpha$ and $q_\beta$ force contradictory information about $\dot{d}\cap x$.
	
	We do so by recursion, iterating over all elements of $\PP$. Enumerate the elements of $\PP$ by $\langle p_\delta \mid \delta < \mu\rangle$. Set $D_{\alpha\beta}^0 = \emptyset$, and for all $\alpha < \nu$, let $q_\alpha^0 = q$.
	First we consider the successor stage $\delta+1$. Applying Claim \ref{cl:singsplitting2} to $p_{\delta}$ and the sequence $\langle q^\delta_\alpha \mid \alpha < \nu\rangle$, we obtain $\langle q^{\delta+1}_\alpha \mid \alpha < \nu\rangle$ with $q^{\delta+1}_\alpha \leq q^\delta_\alpha$ for each $\alpha < \nu$, $x_{\delta+1} \in \mc{P}_{\nu^+}(\la)$, and strengthenings $p_{\alpha\beta}^{\delta} \leq p_{\delta}$ for all $\alpha < \beta < \nu$, so that $(p_{\alpha\beta}^{\delta}, q^{\delta+1}_\alpha)$ and $(p_{\alpha\beta}^{\delta}, q^{\delta+1}_\beta)$ force contradictory information about $\dot{d}\cap x_{\delta+1}$. We set $D_{\alpha\beta}^{\delta+1} = D_{\alpha\beta}^\delta \cup \{p_{\alpha\beta}^{\delta}\}$.
	
	At each limit stage $\gamma$, we take unions for $x_\gamma$ and each $D_{\alpha\beta}^\gamma$, and lower bounds to produce each $q^\gamma_\alpha$. Since $\Q$ is $\mu^+$-closed, we will always be able to take these lower bounds. Let $q_\alpha$ be a lower bound of the sequence $\langle q^\delta_\alpha \mid \delta < \mu\rangle$. Finally, setting $x = \bigcup_{\delta<\mu}x_\delta$, we obtain the desired objects.
	Repeating this construction, we can obtain $(q_\sigma, x_\sigma)$ and $\langle D_\sigma^{\alpha\beta}\mid \alpha < \beta <\nu\rangle$ as desired for any $\sigma \in \nu^{<\tau}$ and any $\alpha < \beta < \nu$.
	
	Let $x = \bigcup_{\sigma \in \nu^{<\tau}} x_\sigma$. For all $f \in \nu^\tau$, by the closure of $\Q$ we can find $q'_f$ such that $q'_f\leq q_{f\uhr \alpha}$ for all $\alpha<\tau$. Let $G$ be generic for $\PP$, and note that it will intersect every $D^{\alpha\beta}_{\sigma}$. Working in $V[G]$, let $q_f \leq q_f'$ be a strengthening that decides $\dot{d}\cap x$.
	
	Let $f \neq g$ be elements of $\nu^\tau$, and let $\alpha$ be the largest such that $f\uhr \alpha = g\uhr \alpha$. Let $\zeta = f(\alpha)$, and $\xi = g(\alpha)$. Since $G$ intersects $D^{\zeta\xi}_{f\uhr \alpha}$, we see that $q_{(f\uhr \alpha) \cat \zeta}$ and $q_{(f\uhr \alpha) \cat \xi}$ force contradictory information about $\dot{d}\cap x_{f\uhr \alpha} \subseteq \dot{d}\cap x$. Since $f$ extends $(f\uhr \alpha) \cat \zeta$ and $g$ extends $(f\uhr n) \cat \xi$, it follows that $q_f$ and $q_g$ disagree on $\dot{d} \cap x$.
	
	We conclude that each $f \in \nu^\tau$ forces a distinct value for $\dot{d}\cap x$, so there must be at least $\nu^\tau > \nu$ potential values for $\dot{d}\cap x$. But since $\dot{d}$ is thinly $\nu^+$-approximated, there are at most $\nu$ potential values for $\dot{d}$ and thus for $\dot{d}\cap x$. This is a contradiction. 
\end{proof}

To obtain the strong tree property and ITP at successors of singular cardinals, as in \cite{NeemanTPNw+1} and \cite{HachtmanITPNw+1} we will use the concept of systems of branches.

\begin{defn}\label{def:system1}
	Let $D \subseteq \Ord, \rho \in \Ord$, and $I$ be an index set. A \emph{system on $D \times \rho$} is a family $\langle R_s \rangle_{s\in I}$ of transitive, reflexive relations on $D \times \rho$ such that:
	\begin{enumerate}
		\item If $(\alpha, \xi) R_s (\beta, \zeta)$ and $(\alpha, \xi) \neq (\beta, \zeta)$, then $\alpha < \beta$.
		\item If $(\alpha_0, \xi_0)$ and $(\alpha_1, \xi_1)$ are both $R_s$-below $(\beta, \zeta)$, then $(\alpha_0, \xi_0)$ and $(\alpha_1, \xi_1)$ are comparable in $R_s$.
		\item For every $\alpha < \beta$ both in $D$, there are $s\in I$ and $\xi, \zeta \in \rho$ so that $(\alpha, \xi) R_s (\beta, \zeta)$.
	\end{enumerate}

	A \emph{branch} through $R_s$ is a subset of $D\times \rho$ that is linearly ordered by $R_s$ and downwards $R_s$-closed. Note that each branch can be viewed as a partial function $b:D\rightharpoonup \rho$. A branch through $R_s$ is \emph{cofinal} if its domain is cofinal in $D$. A \emph{system of branches} through $\langle R_s\rangle_{s\in I}$ is a family $\langle b_\eta\rangle_{\eta \in J}$ so that each $b_\eta$ is a branch through some $R_{s(\eta)}$, and $D = \bigcup_{\eta \in J} \dom(b_\eta)$.
\end{defn}

We have the following branch lemma for systems of branches. This lemma was originally stated for a singular cardinal $\nu$ of countable cofinality, but the same proof applies for any cofinality.

\begin{lemma}\label{lem:1cardbranch}\cite[Lemma 3.3 and Remark 3.4]{NeemanTPNw+1}
	Let $\nu$ be a strong limit singular cardinal, and $\langle R_s\rangle_{s\in I}$ be a system in $V$ on $D\times \rho$, with $D$ cofinal in $\nu^+$. Let $\PP$ be a poset in $V$, and let $\ka < \nu$ be a regular cardinal above $\max(|I|, \rho)^+$, so that:
	\begin{enumerate}
		\item The empty condition in $\PP$ forces that there exists a system $\langle b_\eta\rangle_{\eta \in J}$ of branches through $\langle R_s \rangle_{s\in I}$, with $|J|^+ < \ka$.
		\item For $\la = \max(|I|,|J|,\rho)^+ < \ka$, there is a poset $\PP^\la$ adding $\la$ mutually generic filters for $\PP$ such that $\PP^\la$ is $<\ka$-distributive.
	\end{enumerate}
	Then there exists $\eta$ so that $b_\eta$ is cofinal and belongs to $V$. In particular, there is $s \in I$ so that in $V$, $R_s$ has a cofinal branch.
\end{lemma}

\section{ITP below \texorpdfstring{$\aleph_{\w+1}$}{ℵ\_ω} with the strong tree property at \texorpdfstring{$\aleph_{\w+1}$}{ℵ\_ω+1}}\label{s:TPuptoNw+1}

Using the techniques from the previous section, we show that in the model of \cite{NeemanTPNw+1}, $\ITP$ holds at $\aleph_n$ for all $2 \leq n < \w$, and the strong tree property holds at $\aleph_{\w+1}$.

\subsection{The Forcing}\label{s:theforcing}
First we recall Neeman's construction. See \cite{NeemanTPNw+1} for a more detailed exposition. (In this paper we adopt the convention that a $\ka$-closed forcing is a forcing where every sequence of length less than $\ka$ has a lower bound; \cite{NeemanTPNw+1} uses a different convention that includes sequences of length $\ka$. So in particular, we will refer to $\mu^+$-closed forcings where \cite{NeemanTPNw+1} refers to $\mu$-closed forcings.)

Let $\langle \ka_n \mid 2 \leq n < \w\rangle$ be an increasing sequence of indestructibly supercompact cardinals with supremum $\nu$, and suppose there is a partial function $\phi$ such that for all $2<n<\w$, $\phi\uhr \ka_n$ is an indestructible Laver function for $\ka_n$. That is, for each $A \in V$, $\gamma \in \Ord$, and every extension $V[E]$ of $V$ by a $<\ka_n$-directed closed forcing, there is a $\gamma$-supercompactness embedding $\pi$ in $V[E]$ with critical point $\ka_n$, such that $\pi\uhr \Ord \in V$, $\pi(\phi)(\ka_n)=A$, and the next point in $\dom(\pi(\phi))$ above $\ka_n$ is greater than $\gamma$.
Such a function can be arranged via the standard construction for Laver indestructibility; see \cite[Section 4]{NeemanTPNw+1} for details.
By restricting the domain if necessary, we can assume that for all $\alpha \in \dom(\phi)$, if $\gamma$ is in $\dom(\phi)\cap \alpha$ then $\phi(\gamma) \in V_\alpha$.

The forcing consists of three parts: Cohen forcings to ensure that $2^{\ka_n} = \ka_{n+2}$ for each $n \geq 2$, Laver preparation to ensure that each $\ka_n$ is indestructibly generically supercompact, and Mitchell-style collapses to make each $\ka_n$ become $\aleph_n$. The situation is somewhat complicated by the fact that in order to obtain ITP (or even just the tree property) at $\aleph_{\w+1}$ later on, we cannot fix the cardinal that will become $\aleph_1$ in advance. As in \cite{NeemanTPNw+1}, we define a set Index that gives us the possibilities for the new $\aleph_1$, and define our poset using a lottery sum over the elements of Index.

First we define the Cohen piece of the forcing. 
\begin{defn}
	For $n \geq 2$, let $\A_n \defeq \Add(\ka_n, \ka_{n+2})$. Let $\ka_0$ denote $\w$, and set $\A_0 \defeq \Add(\w, \ka_2)$. Let $\A_1 = \sum_{\mu \in \text{Index}} \Add(\mu^+, \ka_3)$. (The set Index will be defined below, as it relies on the initial stage of both the Cohen and Laver parts of the forcing.) We use $\ka_1$ to refer to the value of $\mu^+$ chosen by a fixed generic. Let $\A$ be the full support product of $\A_n$ for $n < \w$.
\end{defn}
We will use $A_{[n,m]}$ to denote $\prod_{n\leq i\leq m}A_i$; half-open and open intervals are defined similarly. $\A\uhr \alpha$ denotes $\A_{[0,n)}\times \A_n\uhr \alpha$, where $n$ is the least such that $\alpha \leq \ka_{n+2}$.

Next we define the Laver preparation.
\begin{defn}\label{def:laver}
	We define two posets, $\BB$ in $V$ and $\UU$ in the extension of $V$ by $\A$, recursively as follows.
	\begin{enumerate}
		\item Conditions $p$ in $\BB$ are functions with $\dom(p) \subseteq \nu$ so that for every inaccessible cardinal $\alpha$, $|\dom(p)\cap \alpha | < \alpha$. In particular, $|\dom(p)\cap \ka_{n+2}| < \ka_{n+2}$ for each $n$.
		\item If $\alpha \in \dom(p)$, then $\alpha$ is inaccessible, not equal to any $\ka_n$, $\alpha \in \dom(\phi)$, and $\phi(\alpha)$ is an $(\A\uhr \alpha)\ast (\dot{\UU}\uhr \alpha)$ name for a poset forced to be $<\alpha$-directed closed.
		\item $p(\alpha)$ is an $(\A\uhr \alpha ) \ast (\dot{\UU}\uhr \alpha )$ name for a condition in $\phi(\alpha)$.
		\item $p^* \leq p$ in $\BB$ if and only if $\dom(p^*)\supseteq \dom(p)$ and for each $\alpha \in \dom(p)$, $(\emptyset, p^*\uhr \alpha)$ forces in $(\A\uhr \alpha) \ast (\dot{\UU}\uhr\alpha)$ that $p^*(\alpha) \leq p(\alpha)$.
		\item Let $A$ be generic for $\A$. $\UU = \dot{\UU}[A]$ has the same conditions as $\BB$. The order is given by $p^* \leq p$ if and only if $\dom(p^*) \supseteq \dom(p)$ and there exists a condition $a^*$ in $A$ so that for every $\alpha \in \dom(p)$, $(a^*\uhr \alpha, p^*\uhr \alpha)$ forces in $(\A\uhr \alpha) \ast (\dot{\UU}\uhr\alpha)$ that $p^*(\alpha) \leq p(\alpha)$.
	\end{enumerate}
\end{defn}

Let $\UU_0 = U\uhr \ka_2$, and for all $n > 0$ let $\UU_n = \UU\uhr[\ka_{n+1}, \ka_{n+2})$. We define $\UU_{[0,n]} = \UU\uhr\ka_{n+2}$; note that $\UU_{[0,n]}$ is a poset in $V[A_{[0,n]}]$.

\begin{defn}
	Let $\beta < \nu$ and let $F$ be a filter in $\A\ast \dot{\UU}\uhr \beta$. We define $\BB^{+F}\uhr[\beta, \nu)$ to consist of conditions $p\in \BB$ with $\dom(p) \subseteq [\beta, \nu)$ with the ordering $p^* \leq p$ if and only if $\dom(p^*) \supseteq \dom(p)$ and there exists $(a,b) \in F$ so that for every $\alpha \in \dom(p)$, $(a\uhr \alpha, b\cup p^*\uhr \alpha)$ forces $p^*(\alpha) \leq p(\alpha)$.
\end{defn}

\begin{rem}
	Given a generic $A \ast U$ for $\A \ast \dot{\UU}$ over $V$, with $F = A\uhr \beta \ast U\uhr \beta$, $U \uhr [\beta,\nu)$ may not be a generic filter for $\BB^{+F} \uhr [\beta,\nu)$. However, we can force to add a filter $G \subseteq U\uhr [\beta,\nu) $ that is generic for $\BB^{+F}\uhr [\beta,\nu)$. We call this forcing the \emph{factor forcing} refining $U\uhr \beta$ to a generic for $\BB^{+F}\uhr[\beta,\nu)$.
\end{rem}

\begin{fact}\cite[Claim 4.7]{NeemanTPNw+1}\label{fact:Brefinementisdirclosed}
	Let $\bar{\beta} \leq \beta$. Suppose that $\bar{F}$ is generic for $\A\uhr\bar{\beta}\ast\dot{\UU}\uhr\bar{\beta}$ over $V$. Then $\BB^{+\bar{F}}\uhr [\beta,\nu)$ is $<\beta$ directed closed in $V[\bar{F}]$.
\end{fact}

\begin{fact}\cite[Lemma 4.12]{NeemanTPNw+1}\label{fact:genindsc}
	Let $n < \w$ and let $A\ast U_{[0,n]}$ be generic for $\A\ast \dot{\UU}_{[0,n]}$ over $V$. Then in $V[A][U_{[0,n]}]$, $\ka_{n+2}$ is generically supercompact. Moreover, this supercompactness is indestructible under forcing with $<\ka_{n+2}$-directed-closed forcings in $V[A_{[0,n]}][U_{[0,n]}]$. The forcing poset producing this embedding is $\Add(\ka_n, \pi(\ka_{n+2}))^V\times \Add(\ka_{n+1}, \pi(\ka_{n+3}))^V$.
	
	By this we mean the following. Let $\hat{A}$ be generic for $\Add(\ka_n, \pi(\ka_{n+2}))^V\times \Add(\ka_{n+1}, \pi(\ka_{n+3}))^V$. Let $\la \geq \ka_{n+2}$, and let $G$ be generic for some $\ka_{n+2}$-directed-closed forcing in $V[A_{[0,n]}][U_{[0,n]}]$. Then in the generic extension $V[A][U_{[0,n]}][G][\hat{A}]$ there is an elementary embedding $\pi: V[A][U_{[0,n]}][G] \to V^*[A^*][U^*_{[0,n]}][G^*]$ with critical point $\ka_{n+2}$, such that $\pi(\ka_{n+2}) > \la$, $\pi\uhr\Ord \in V$, and $V^*[A^*][U^*_{[0,n]}][G^*]$ is closed under $\la$-sequences in this generic extension.
	
\end{fact}

\begin{defn}
	We define the set Index as the set of all $\mu < \ka_2$ so that:
	\begin{enumerate}
		\item $\mu$ is a strong limit cardinal of cofinality $\w$ and $\dom(\phi)$ has a largest point $\la$ below $\mu$.
		\item Over any extension $V[E]$ of $V$ by a $\mu^+$ closed poset, the further extension by $\A_0\uhr \la \ast \dot{U}_0\uhr \la + 1$ does not collapse $(\mu^+)^V$.
		\item $\A_0 \uhr \la \ast \dot{U}_0\uhr \la +1$ has size at most $\mu^+$.
	\end{enumerate}
\end{defn}

We note that there are many elements in $\Index$. In particular, consider the case where $\gamma$ is a strong limit cardinal with cofinality $\omega$ above $\ka_2$ and $i:V \to V^*$ is some $\gamma$-supercompactness embedding with critical point $\ka_2$. Then the largest point of $\dom(i(\phi))$ below $\gamma$ will be $\ka_2$; by elementarity there are many $\mu$ satisfying the first condition. If $|\phi(\la)| < \mu$, then the second and third conditions will be satisfied, since $\A_0\uhr \la\ast\dot{U}\uhr \la+1$ will have size less than $\mu$ and so will not be able to collapse $\mu^+.$

Finally, we define the collapses.
\begin{defn}\label{def:collapses}
	For each $n < \w$, define $\CC_n$ in $V$ as follows. Conditions in $\CC_n$ are functions $p$ such that:
	\begin{enumerate}
		\item $\dom(p)$ is contained in the interval $(\ka_{n+1}, \ka_{n+2})$ and $|\dom(p)| < \ka_{n+1}$.
		\item For each $\alpha \in \dom(p)$, $p(\alpha)$ is an $(\A\uhr \alpha)\ast (\dot{U}\uhr \ka_{n+1})$ name for a condition in the poset $\Add(\ka_{n+1}, 1)$ of the extension by $(\A\uhr \alpha)\ast (\dot{U}\uhr \ka_{n+1})$.
	\end{enumerate}

	We define the ordering by $p^* \leq p$ if and only if $\dom(p^*) \supseteq \dom(p)$ and for each $\alpha \in \dom(p)$, the empty condition in $(\A\uhr \alpha)\ast (\dot{U}\uhr \ka_{n+1})$ forces that $p^*(\alpha) \leq p(\alpha)$.
	
	Let $\CC$ be the full support product of each $\CC_n$.
\end{defn}

Note that $\CC_0$ is defined using $\ka_1$, which is not a priori known. When we have access to a generic for $\A_1$, we will use the $\ka_1$ determined by this generic; if not, we will view it as a parameter in the definition. We will occasionally refer to this poset as $\CC_0(\ka_1)$ if the choice of $\ka_1$ is not clear from context.

\begin{defn}
	For a filter $F \subseteq \A\ast\dot{U}$ define the \emph{enrichment of $\CC$ to $F$}, denoted by $\CC^{+F}$, to be the poset with the same conditions as $\CC$ but the order given by $p^* \leq p$ if and only if there exists a condition $(a, u) \in F$ so that for each $\alpha \in \dom(p)$, $(a\uhr \alpha, u\uhr \ka_i) \forces_{\A\uhr \alpha\ast \dot{U}\uhr \ka_{i}} p^*(\alpha) \leq p(\alpha)$, where $i$ is the largest such that $\ka_i \leq \alpha$.
\end{defn}

\begin{rem}
	Let $A \ast U$ be generic for $\A\ast \UU$, $F = A\uhr \ka_n \ast U\uhr \ka_n$, and let $S$ be generic for $\CC^{+A\ast U}$. As before, we can find a factor forcing refining $S\uhr[\ka_n, \nu)$ to a generic for $\CC^{+F}\uhr[\ka_n, \nu)$.
\end{rem}

\begin{fact}\cite[Claim 4.15]{NeemanTPNw+1}\label{fact:Crefinementdirclosed}
	\begin{enumerate}
		\item Let $F$ be generic for $\A\uhr\beta\ast\dot{\UU}\uhr\beta$ for $\beta \leq \ka_{n+1}$. Then $\CC^{+F}\uhr[\ka_{n+1},\nu)$ is $<\ka_{n+1}$ directed closed in $V[F]$.
		\item Let $\beta \in (\ka_{n+1}, \ka_{n+2})$ and let $F$ be generic for $\A\uhr\beta\ast\dot{U}\uhr\ka_{n+1}.$ Then $\CC^{+F}\uhr[\beta,\ka_{n+2})$ is $<\ka_{n+1}$ directed closed in $V[F]$.
	\end{enumerate}
\end{fact}

Let $A$ be generic for $\A$, $U$ be generic for $\UU$ over $V[A]$, $S$ be generic for $\CC^{+A\ast U}$ over $V[A][U]$, and $e$ be generic for $\Coll(\w, \mu)$ over $V[A][U][S]$, noting that $\mu$ is determined by the generic for $\A_1$ included in $A$. The final model is $V[A][U][S][e]$.

The cardinal structure in this model is what one would expect.
\begin{fact}
	In $V[A][U][S][e]$, the following properties hold.
	\begin{itemize}
		\item $\ka_n = \aleph_n$ for each $n$.
		\item $2^{\ka_n} = \ka_{n+2}$ for each $n$.
		\item $2^\w = \ka_2$.
	\end{itemize}
\end{fact}
\begin{fact}\cite[Lemma 4.24]{NeemanTPNw+1}\label{fact:covprop}
	$V$ has the $<\ka_n$-covering property with respect to $V[A][U][S][e]$ for all $n \geq 2$.
\end{fact}

\subsection{The strong tree property at the successor of a singular}

To obtain the strong tree property at $\aleph_{\w+1}$, we will generalize \cite[Lemma 3.10]{NeemanTPNw+1} to the strong tree property. For full generality, and to allow us to use it in later results, we will prove the lemma for the successor of a singular cardinal of any cofinality. The main obstacle here is that the strong tree property is a global property that must hold for every $\la$; we solve this by bringing in an auxiliary poset that collapses $\la$ to $\nu^+$.

\begin{thm}\label{thm:TP2cardfullNA}
	Let $\tau$ be a regular cardinal. Let $\langle \ka_{\rho} \mid \rho < \tau\rangle$ be a continuous sequence of regular cardinals above $\tau$, with supremum $\nu$. Let $\Index \subseteq \ka_0$, and fix $\rho' < \tau$. For each $\mu\in \Index$, let $\LL_\mu$ be a forcing poset
	of size $\leq \ka_{\rho'}$.
	Let $R$ be a rank-initial segment of $V$ satisfying a large enough fragment of ZFC. For all $\la \geq \nu^+$ with $\la^\nu = \la$, let $\bbK = \bbK_\la$ be a poset preserving all cardinals $\leq \nu^+$ that forces $|\la|=\nu^+$.
	Let $K$ be generic for $\bbK$. Suppose that:
	\begin{enumerate}
		\item For each $X \prec R$ with $|X| = \la$ and $\la \subseteq X$ such that $\bbK$ is in the transitive part of $X$, let $\bar{V} = \bar{V}_X$ be the transitive collapse of $X$. For stationarily many $X$, in $V[K]$ there exists a $\nu^+$-Knaster poset $\HH = \HH_X$ such that $\HH$ forces the existence of a generic $\nu^+$-supercompactness embedding $i: \bar{V}[K] \to \bar{M}[K^*]$ with critical point $\ka_0$ such that $\nu \in i(\Index)$, and a set $L$ such that $L$ is generic over $\bar{M}[K^*]$ for $i(\LL)(\nu)$.
		
		\item In $V[K]$, for all ordinals $\rho < \tau$, there is a generic $\nu^+$-supercompactness embedding $j_{\rho+2}$ with domain $V[K]$ and critical point $\ka_{\rho+2}$, added by a poset $\FF$ such that the full support power $\FF^{\ka_\rho}$ is $<\ka_{\rho}$-distributive in $V[K]$.
	\end{enumerate}
	Then there is $\mu \in \Index$ such that the strong tree property holds at $\nu^+$ in the extension of $V$ by $\LL_\mu$.
\end{thm}
\begin{proof}
	Suppose not. Then for every $\mu \in \Index$, there exists $\la \geq \nu^+$ regular such that there is a thin $\mc{P}_{\nu^+}(\la)$-list with no cofinal branch in the extension of $V$ by $\LL_\mu$. Applying Fact \ref{fact:TPcandecreaselambda}, by taking an upper bound, we can assume that $\la$ is the same for each $\mu$. Similarly, by increasing $\la$ if necessary we can assume that $\la^\nu = \la$. For each $\mu \in \Index$, let $\dot{d}^\mu$ be a $\LL_\mu$-name for a thin $\mc{P}_{\nu^+}(\la)$-list forced by the empty condition to have no cofinal branch. Assume that the empty condition in $\LL_\mu$ forces the $z$-th level of $\dot{d}^\mu$ to be indexed by $\langle \dot{\sigma}_z^\mu (\zeta) \mid \zeta < \nu\rangle$. We work in $V[K]$. Note that $\la$ is collapsed, but we can still examine cofinal branches through a thin $\mc{P}_{\nu^+}(\la)$-list when $\la$ is merely an ordinal. In fact, since $\la$ has cardinality and cofinality $\nu^+$, $\mc{P}_{\nu^+}(\la)$ is order isomorphic to $\mc{P}_{\nu^+}(\nu^+)$. We can thus assume that, in $V[K]$, $\dot{d}^\mu$ is a $\LL_\mu$-name for a thin $\mc{P}_{\nu^+}(\nu^+)$-list. Since $\nu^+$ is cofinal (actually club) in $\mc{P}_{\nu^+}(\nu^+)$, we can further assume that $\dot{d}^\mu$ is indexed by ordinals.
	
	Let $I = \{(r, \mu) \mid \mu \in \Index, r \in \LL_\mu\}$. For $s = (r,\mu)$, let $R_s$ be the relation $(\alpha, \zeta) R_s (\beta, \xi)$ iff $\alpha \leq \beta$ and $r \forces_{\LL_{\mu}} \dot{\sigma}^\mu_\alpha(\zeta) = \dot{\sigma}^\mu_{\beta}(\xi)\cap \alpha$.
	
	Let $X, \bar{V}$, and $\HH$ be as in condition (1), with the function $\mu \mapsto \dot{d}^\mu$ contained in $X$. Let $G$ be generic for $\HH$ over $V[K]$, and let $i, L \in V[K][G]$ be as in condition (1). Noting that $\nu \in i(\Index)$, we can apply the image of this map to $\nu$ to obtain a $i(\LL)(\nu)$-name for a thin $\mc{P}_{\nu^+}(\la)$-list, that we will denote by $i\dot{d}^\nu$. Let $id^\nu = (i\dot{d}^{\nu})_L \in V[K][G]$, and denote the $\alpha$-th level with $i\sigma^\nu_\alpha$ similarly.
	
	Let $\gamma \in i(\nu^+)$ such that $\gamma > i''\nu^+$. For each $\alpha < \nu^+$, there exists $\rho_\alpha$ with $\rho' \leq \rho_\alpha < \tau$ and $\xi_\alpha < i(\ka_{\rho_\alpha})$ such that $i\sigma_{i(\alpha)}^\nu(\xi_\alpha) = i\sigma_\gamma^\nu(0)\cap i(\alpha)$.
	Let $\dot{\rho}_\alpha$ be an $\HH$-name for $\rho_\alpha$. Note that the construction of this name does not depend on the choice of generic $G$, since the existence of such a $\rho_\alpha$ is forced by the empty condition; the generic is used for notational convenience.
	
	Working now in $V[K]$, for each $\alpha < \nu^+$, let $h_\alpha \in \HH$ decide the value of $\dot{\rho}_\alpha$ to be some ordinal $\rho_\alpha < \tau$. Since $\tau < \nu^+$ are both regular, and $\HH$ is $\nu^+$-Knaster, in $V[K]$ we can obtain an unbounded set $S \subseteq \nu^+$ and a fixed $\rho < \tau$ such that for any $\alpha,\beta$ both in $S$, $h_\alpha$ and $h_{\beta}$ are compatible, and $\rho_\alpha = \rho_{\beta} = \rho$. Note that by construction, $\rho' \leq \rho$, so $|\LL_\mu| \leq \ka_\rho$.	
	\begin{lemma}
		$\langle R_s \rangle_{s \in I}$ is a system on $S \times \ka_\rho$.
	\end{lemma}
	\begin{proof}
		The first two conditions are trivial, so it suffices to verify the third. Let $\alpha < \beta \in S$. Let $p$ be a common extension of $h_\alpha$ and $h_\beta$, and let $G$ be a generic for $\HH$ over $V[K]$ containing $p$. 
		Work in $V[K][G]$.
		Then there exist $\zeta, \xi < i(\ka_\rho)$ such that $i\sigma^\nu_{i(\alpha)}(\zeta) = i\sigma^\nu_\gamma(0) \cap i(\alpha)$ and $i\sigma^\nu_{i(\beta)}(\xi) = i\sigma^\nu_\gamma(0) \cap i(\beta)$. We conclude that $\bar{M}[K^*] \models \exists \mu \in i(\Index), r\in i\LL_{\mu}$, and $\zeta, \xi < i(\ka_\rho)$ such that $r \forces \sigma^\nu_{i(\alpha)}(\zeta) = \sigma^\nu_{i(\beta)}(\xi)\cap i(\alpha)$. By elementarity, there exists $\mu \in \Index, r\in \LL_{\mu}$, and $\zeta, \xi < \ka_\rho$ such that $r \forces \sigma^\mu_\alpha(\zeta)= \sigma^\mu_{\beta}(\xi)\cap \alpha$. In particular, letting $s = (r,\mu)$, we have that $(\alpha, \zeta) R_s (\beta, \xi)$.
	\end{proof}
	
	\begin{lemma}\label{lem:TP2cardgettingabranch}
		There exists $s \in I$ such that $R_s$ has a cofinal branch.
	\end{lemma}
	\begin{proof}
		By assumption, we have a generic $\nu^+$-supercompactness embedding $j$ with domain $V[K]$ and critical point $\ka_{\rho+3}$, added by a poset $\FF$ with generic $F$ whose full support power is $<\ka_{\rho+1}$-distributive in $V$. In particular, this power is $<\ka_{\rho+1}$-distributive in $V[K]$, so it satisfies hypothesis (2) of Lemma \ref{lem:1cardbranch}. Recall that $|\LL_\mu| < \ka_0$. We can assume that the underlying set of $\LL_\mu$ is an ordinal below the critical point of $j$, so that $j$ fixes $\LL_\mu$ pointwise.
		
		Work in $V[K][F]$. Let $\gamma' \in j(S)$ such that $\gamma' > \sup j''\nu^+$. For each $\delta < \ka_\rho$ and $s = (\mu, p) \in I,$ we define
		\[b_{s,\delta} = \{(\alpha, \zeta) \mid \alpha \in S, \zeta < \ka_\rho, p \forces_{\LL_\mu} j(\dot{\sigma}_\alpha^\mu(\zeta)) = j\dot{\sigma}_{\gamma'}^\mu(\delta) \cap j(\alpha)\}.\]
		
		We claim that $\langle b_{s,\delta} \mid s\in I, \delta < \ka_\rho\rangle$ is a system of branches through $\langle R_s\rangle_{s\in I}$. By construction, each $b_{s,\delta}$ is linearly ordered and downwards closed. It remains to verify that $\bigcup \dom(b_s,\delta) = S$.
		Since $\ka_\rho < \crit(j)$, and $\langle R_s\rangle_{s\in I}$ is a system, we can apply elementarity to conclude that for all $\alpha \in S$ there exists $\zeta, \delta < \ka_\rho$ and $s = (\mu, p) \in j(I) = I$ such that
		\[p \forces_{\LL_{\mu}} j(\dot{\sigma}_\alpha^{\mu}(\zeta)) = j\dot{\sigma}_{\gamma'}^{\mu}(\delta) \cap j(\alpha).\]
		In particular, $\bigcup \dom(b_{s,\delta}) = S$.
		
		Note that this system of branches is in $V[K][F]$. Noting that $|I| \leq \ka_\rho$, we apply Lemma \ref{lem:1cardbranch}, we conclude that there is some $(s, \delta) \in I \times \ka_{\rho}$ such that $b_{s,\delta}$ is cofinal and belongs to $V[K]$.
	\end{proof}
	
	 Now we finish the proof. Let $b_{s,\delta}$ be the cofinal branch from the previous lemma, with $s = (\mu, p)$, and let $L$ be generic for $\LL_\mu$ containing $p$. We define an $\LL_\mu$-name $\dot{\pi}_{s,\delta} = \bigcup\{\dot{\sigma}_\alpha^\mu(\zeta) \mid (\alpha, \zeta) \in b_{s,\delta}\}$. Since $b_{s,\delta}$ is a cofinal branch through $R_s$, $\dot{\pi}_{s,\delta}$ will be a cofinal branch through $d^\mu$ in $V[K][L]$.
	
	Since $\bbK$ is $\ka_{\rho'}$-closed in $V$, and $\LL_\mu$ is $\ka_{\rho'}$-cc, by Lemma \ref{lem:M-S} we see that $\bbK$ is forced by $\LL_\mu$ to have the thin $\nu^+$-approximation property. It follows that $\bbK$ could not have added a cofinal branch to $V[L]$, so $d^\mu$ must have a cofinal branch in $V[L]$, contradicting our original assumption.
\end{proof}

\subsection{ITP below \texorpdfstring{$\aleph_{\w}$}{ℵ\_ω}}\label{section:ITPuptoNw}

Next, we show that in this model, we have ITP at $\aleph_n$ for all $n \geq 2$. The argument follows that of \cite[Section 4]{NeemanTPNw+1}, using the stronger branch preservation lemmas described in Section \ref{section:lemmas}.

\begin{theorem}\label{thm:ITPatNn}
	In the model $V[A][U][S][e]$, $\ITP$ holds at $\aleph_{n+2}$ for all $n$.
\end{theorem}
\begin{proof}
	Let $d = \langle d_z \mid z \in (\mathcal{P}_{\ka_{n+2}}(\la))^{V[A][U][S][e]} \rangle$ be a thin $\mc{P}_{\ka_{n+2}}(\la)$ list in $V[A][U][S][e]$. We wish to show that it has an ineffable branch.

	We repeat the argument of \cite{NeemanTPNw+1}, examining in particular the proof of \cite[Theorem 4.29]{NeemanTPNw+1}. This construction will produce a lifted $\la$-supercompactness embedding with critical point $\ka_{n+2}$ contained in an extension of $V[A][U][S][e]$. We record the definition and some properties of the forcings used to produce this embedding. A detailed exposition can be found in \cite[Section 4]{NeemanTPNw+1}.
	
	Let $F = A\uhr\ka_{n+2} \ast U\uhr \ka_{n+2}$. Let $\PP_1$ be the factor forcing refining $U\uhr[\ka_{n+2}, \nu)$ to a generic $G_1$ for $\BB^{+F}\uhr[\ka_{n+2}, \nu)$, and let $\PP_2$ be the factor forcing refining $S\uhr[\ka_{n+2}, \nu)$ to a generic $G_2$ for $\CC^{+F}\uhr [\ka_{n+2}, \nu)$.
	
	\begin{prop}
		There exists a generic $\la$-supercompactness embedding \[\pi : V[A][U_{[0,n]}][G_1][G_2] \to V^*[A^*][U^*_{[0,n]}][G_1^*][G_2^*]\]
		with critical point $\ka_{n+2}$, contained in the extension of $V[A][U_{[0,n]}][G_1][G_2]$ by the product of $\Add(\ka_n, \pi(\ka_{n+2}))^V$ and $\Add(\ka_{n+1}, \pi(\ka_{n+3}))^V$. 
	\end{prop}
	\begin{proof}
		$\BB^{+F}\uhr[\ka_{n+2},\nu)$ and $\CC^{+F}\uhr[\ka_{n+2},\nu)$ are in $V[F]$, and by Fact \ref{fact:Brefinementisdirclosed} and Fact \ref{fact:Crefinementdirclosed} they are $<\ka_{n+2}$ directed closed in $V[F]$. So using Fact \ref{fact:genindsc}, we obtain the desired embedding.
	\end{proof}
	Consider the map $\beta \mapsto \CC^{+A\uhr\beta\ast U\uhr\ka_{n+1}}\uhr [\beta, \ka_{n+2})$ described in the second part of Fact \ref{fact:Crefinementdirclosed}, defined for $\beta \in (\ka_{n+1}, \ka_{n+2})$. Note that the image of $\beta$ under this map is in $V[A\uhr\beta][ U\uhr\ka_{n+1}].$ Consider the image of this map under $\pi$; this takes each $\beta \in (\ka_{n+1}, \pi(\ka_{n+2}))$ to a forcing contained in $V^*[A^*\uhr \beta][ U^*\uhr \ka_{n+1}]$.
	Let $\PP_3$ denote the image of $\ka_{n+2}$ under this new map; i.e., $\PP_3 = \pi(\CC)^{+A\uhr\ka_{n+2}\ast U\uhr\ka_{n+1}}\uhr [\ka_{n+2}, \pi(\ka_{n+2}))$. Note that $\PP_3 \in V^*[A^*\uhr\ka_{n+2}][U^*\uhr\ka_{n+1}] = V^*[A\uhr\ka_{n+2}][U\uhr\ka_{n+1}]$. Since $V^* \subseteq V[A_{[n+2, \omega)}]$ (see the proof of \cite[Lemma 4.12]{NeemanTPNw+1}), we conclude that $\PP_3 \in V[A_{[n+2, \omega)}][A\uhr\ka_{n+2}][U\uhr\ka_{n+1}]$, so in particular $\PP_3\in V[A][U][S][e]$. Let $G_3$ be generic for $\PP_3$ over $V[A][U][S][e]$.
	
	\begin{prop}
		The embedding $\pi$ restricts to an embedding from $V[A][U][S_{[n+1,\w)}]$ to $V^*[A^*][U^*][S^*_{[n+1,\w)}]$; this restriction extends to an embedding from $V[A][U][S][e]$ to $V^*[A^*][U^*][S^*][e]$ (which we will also denote by $\pi$). This embedding is contained in the generic extension of $V[A][U][S][e]$ by the product of $\PP_1$, $\PP_2$, $\Add(\ka_n, \pi(\ka_{n+2}))^V$, $\Add(\ka_{n+1}, \pi(\ka_{n+3}))^V$, and $\PP_3$.
	\end{prop}
	\begin{proof}
		Let $U^*_{[n+1,\w)}$ be the upwards closure of $G^*_1$ in $\pi(\UU_{[n+1,\w)})$. It is generic for $\pi(\UU_{[n+1,\w)})$ over $V^*[A^*][U^*][S^*_{[n+1,\w)}]$. Let $U^*$ be the concatenation of $U^*_{[0,n]}$ and $U^*_{[n+1,\w)}$. Then $U^*$ is generic for $\pi(\UU)$ over $V^*[A^*]$, and $G^*_2$ is generic over $V[A^*][U^*]$. We can restrict $\pi$ to an embedding (in a mild abuse of notation, we will denote all restrictions and extensions by $\pi$) from $V[A][U][G_2]$ to $V^*[A^*][U^*][G_2^*]$.
		
		Let $S^*_{[n+1,\w)}$ be the upwards closure of $G_2^*$ in $\pi(\CC)^{+A^*\ast U^*}_{[n+1,\w)}$. As before we can restrict $\pi$ further to an embedding from $V[A][U][S_{[n+1,\w)}]$ to $V^*[A^*][U^*][S^*_{[n+1,\w)}]$. Noting that $e$ and $S_{[0,n-1]}$ are generic for forcings of size less than $\ka_{n+2}$, and are thus below the critical point of $\pi$, without any further forcing we can extend $\pi$ to an embedding from $V[A][U][S_{[n+1,\w)}][S_{[0,n-1]}][e]$ to $V^*[A^*][U^*][S^*_{[n+1,\w)}][S_{[0,n-1]}][e].$
		
		Finally, let $G_3^+$ be the upwards closure of $G_3$ in $\pi(\CC_n)^{+A^*\ast U^*}\uhr[\ka_{n+2}, \pi(\ka_{n+2}))$. Let $S^*_n = S_n \times G_3^+$. This is generic for $\pi(\CC_n^{+A\ast U}) = \pi(\CC_n)^{+A^*\ast U^*}$. We can extend $\pi$ to an embedding from $V[A][U][S][e] = V[A][U][S_{[n+1,\w)}][S_{[0,n-1]}][e][S_n]$ to $V[A^*][U^*][S^*_{[n+1,\w)}][S_{[0,n-1]}][e][S^*_n] = V^*[A^*][U^*][S^*][e].$
	\end{proof}
	Let $\hat{A}_{n}$ be a generic filter for $\Add(\ka_n, \pi(\ka_{n+2}))^V$ and let $\hat{A}_{n+1}$ be generic for $\Add(\ka_{n+1}, \pi(\ka_{n+3}))^V$. Since all posets are contained in $V[A][U][S][e]$, we can add them in any order. We will add $\hat{A}_{n+1}$ first, followed by $G_1 \times G_2 \times G_3$ and then $\hat{A}_n$.
	
	Since the posets adding the embedding are not defined over the full model, we need some facts about the construction to determine their properties.
	
	\begin{fact}\cite[Lemma 4.26]{NeemanTPNw+1}\label{fact:laterforcingsdist}
		For all $n<\w$, all sequences of ordinals of length $<\ka_{n+1}$ in $V[A][U][S][e]$ belong to $V[A\uhr \ka_{n+2}][U\uhr \ka_{n+1}][S\uhr \ka_{n+1}][e]$.
	\end{fact}

	\begin{fact}\cite[Remark 4.27]{NeemanTPNw+1}\label{fact:extraforcingsdist}
		If a forcing $\Q$ is $\ka_{n+1}$-closed in $V$, then it is $<\ka_{n+1}$-distributive over $V[A][U][S][e]$.
	\end{fact}
	
	With these facts in hand, we record some properties of these posets.
	\begin{prop}\label{fact:p1p2p3closed}\cite[Claims 4.30 and 4.31]{NeemanTPNw+1}
		$\PP_1$, $\PP_2$, and $\PP_3$ are $\ka_{n+1}$ closed in $V[A][U][S\uhr [\ka_{n+1},\nu)][\hat{A}_{n+1}]$.
	\end{prop}
	\begin{proof}
		First we show that each poset is closed in $V[A][U][S\uhr [\ka_{n+1},\nu)]$.
		
		By Fact \ref{fact:Brefinementisdirclosed}, every decreasing sequence of $\PP_1$ with length $<\ka_{n+2}$ in $V[F]$ has a lower bound. By Fact \ref{fact:laterforcingsdist}, every decreasing sequence of length $<\ka_{n+1}$ belonging to $V[A][U][S\uhr[\ka_{n+1}, \nu)]$ must be in $V[A\uhr \ka_{n+2}][U\uhr \ka_{n+1}][S\uhr\ka_{n+1}]$. The sequence came from $V[A][U][S\uhr[\ka_{n+1}, \nu)]$, and that model and $V[A\uhr \ka_{n+2}][U\uhr \ka_{n+1}][S\uhr\ka_{n+1}]$ are mutually generic extensions of $V[A\uhr \ka_{n+2}][U\uhr \ka_{n+1}]$. Since the sequence is present in both models, it must be be in $V[A\uhr \ka_{n+2}][U\uhr \ka_{n+1}] \subseteq V[F],$ and thus have a lower bound in $\PP_1$. We conclude that $\PP_1$ is $\ka_{n+1}$-closed in $V[A][U][S\uhr [\ka_{n+1},\nu)]$ as desired.
		
		The proof for $\PP_2$ is analogous, using Fact \ref{fact:Crefinementdirclosed} in place of Fact \ref{fact:Brefinementisdirclosed}.
		
		Finally we verify that $\PP_3$ is closed. By the second part of Fact \ref{fact:Crefinementdirclosed} applied in $V^*[A^*][U^*]$, $\PP_3$ is $\ka_{n+1}$-closed in $V^*[A^*\uhr\ka_{n+2}][U^*\uhr\ka_{n+1}] = V^*[A\uhr\ka_{n+2}][U\uhr\ka_{n+1}]$. $V^*$ is $\ka_{n+2}$-closed in $V[A_{[n+2,\w)}]$, so $V^*[A\uhr\ka_{n+2}][U\uhr\ka_{n+1}]$ is $\ka_{n+2}$-closed in $V[A_{[n+2,\w)}][A\uhr\ka_{n+2}][U\uhr\ka_{n+1}]$. We conclude that $\PP_3$ is $\ka_{n+1}$-closed in $V[A_{[n+2,\w)}][A\uhr\ka_{n+2}][U\uhr\ka_{n+1}]$.
		By Fact \ref{fact:laterforcingsdist}, any sequence of ordinals of length $<\ka_{n+1}$ in $V[A][U][S\uhr[\ka_{n+1},\nu)]$ is in $V[A\uhr\ka_{n+2}][U\uhr\ka_{n+1}]$, and thus in $V[A_{[n+2,\w)}][A\uhr\ka_{n+2}][U\uhr\ka_{n+1}]$. So any descending sequence of conditions of $\PP_3$ that is in $V[A][U][S\uhr[\ka_{n+1},\nu)]]$ is in $V[A_{[n+2,\w)}][A\uhr\ka_{n+2}][U\uhr\ka_{n+2}]$, and thus has a lower bound. We conclude that $\PP_3$ is $\ka_{n+1}$-closed in $V[A][U][S\uhr [\ka_{n+1},\nu)]$.
		
		Since $\hat{A}_{n+1}$ is $\ka_{n+1}$-closed in $V$, by Fact \ref{fact:extraforcingsdist} it is $\ka_{n+1}$-distributive over $V[A][U][S\uhr [\ka_{n+1},\nu)]$. It follows that $\PP_1$, $\PP_2$, and $\PP_3$ remain $\ka_{n+1}$-closed in $V[A][U][S\uhr [\ka_{n+1},\nu)][\hat{A}_{n+1}]$.
	\end{proof}

	\begin{prop}\label{fact:ccextension}
		$V[A][U][S][e][\hat{A}_{n+1}]$ is a $\ka_{n+1}$-cc extension of the model $V[A][U][S\uhr [\ka_{n+1},\nu)][\hat{A}_{n+1}]$.
	\end{prop}
	\begin{proof}
		The extending poset is $\CC_{[0,n)}^{+A\ast U}\times\Coll(\w,\mu)$. Every piece of this poset has size $<\ka_{n+1}$ except $\CC_{n-1}^{+A\ast U}$, so it is enough to verify that this poset is $\ka_{n+1}$-cc. By Fact \ref{fact:covprop}, $V$ has the $<\ka_{n+1}$-covering property in $V[A][U][S][e]$; since $\hat{A}_{n+1}$ is $<\ka_{n+1}$-distributive, $V$ likewise has the $<\ka_{n+1}$-covering property in $V[A][U][S][e][\hat{A}_{n+1}]$, and thus also in $V[A][U][S\uhr [\ka_{n+1},\nu)][\hat{A}_{n+1}]$. A standard $\Delta$-system argument, using Lemma \ref{lem:covtocc} and the fact that $\ka_{n+1}$ is inaccessible in $V[A][U][S\uhr [\ka_{n+1},\nu)][\hat{A}_{n+1}]$, yields that the poset has the $\ka_{n+1}$ chain condition.
	\end{proof}
	
	By Lemma \ref{lem:embeddingtobranch}, $d$ has a cofinal branch $b$ in $V[A][U][S][e][\hat{A}_{n+1}][G_1\times G_2\times G_3][\hat{A}_n]$, that is ineffable if it is contained in $V[A][U][S][e]$. It remains to show that this branch is actually present in $V[A][U][S][e]$. To do this, we will show that at each stage, the forcing has the appropriate approximation property, and then verify that $\dot{b}$ is approximated by that forcing.
	
	\begin{claim}
		The branch $b$ is in $V[A][U][S][e][\hat{A}_{n+1}][G_1\times G_2\times G_3].$
	\end{claim}
	\begin{proof}
		First we show that $\Add(\ka_n, \pi(\ka_{n+2}))^V$ has the $\ka_{n+1}$-approximation property. Note that $\ka_{n+2}$ is collapsed in the extension by $G_1\times G_2 \times G_3$, since $G_3$ will collapse all cardinals in the interval $[\ka_{n+2}, \pi(\ka_{n+2}))$, but the cofinality of $\ka_{n+2}$ in the extension is at least $\ka_{n+1}$.
		
		Since $V[A][U][S][e][\hat{A}_{n+1}]$ is the extension of $V[A][U][S\uhr [\ka_{n+1},\nu)][\hat{A}_{n+1}]$ by a $\ka_{n+1}$-cc forcing, by Lemma \ref{lem:ccpreservesdistributive} we see that $G_1\times G_2 \times G_3$ doesn't add any new sequences of ordinals of length $\ka_{n+1}$ to $V[A][U][S][e][\hat{A}_{n+1}]$. If $n \geq 1$, by Fact \ref{fact:covprop}, $V$ has the $\ka_{n+1}$-covering property in $V[A][U][S][e][\hat{A}_{n+1}][G_1\times G_2\times G_3].$ Applying Lemma \ref{lem:covtocc}, we see that that $\Add(\ka_n, \pi(\ka_{n+2})\cdot 2)^V$, which is the square of $\Add(\ka_n, \pi(\ka_{n+2}))^V$, is $\ka_{n+1}$-cc over this model. Lemma \ref{lem:ccbranch} gives us that the poset has the $\ka_{n+1}$-approximation property. For $n=0$, we see that the square of $\Add(\ka_n, \pi(\ka_{n+2}))^V$ is $\ka_1$-cc in any model where $\ka_1$ is a cardinal.
		
		Let $\dot{b}$ be a name (in $V[A][U][S][e]$) for the branch $b$. We will verify that $\dot{b}$ is $\ka_{n+1}$-approximated by $\Add(\ka_n, \pi(\ka_{n+2}))^V$ over $V[A][U][S][e][\hat{A}_n][G_1\times G_2 \times G_3]$.
		Let $W$ denote $V[A][U][S\uhr[\ka_{n+1},\nu)][\hat{A}_{n+1}]$.
		We note that the product $\PP_1 \times \PP_2 \times \PP_3$ is $\ka_{n+1}$-closed in $W$ by Fact \ref{fact:p1p2p3closed}, and $V[A][U][S][e][\hat{A}_{n+1}]$ is a $\ka_{n+1}$-cc extension of $W$ by Fact \ref{fact:ccextension}. Then by Lemma \ref{lem:ccpreservesdistributive}, forcing with $\PP_1\times\PP_2\times\PP_3$ over $V[A][U][S][e][\hat{A}_{n+1}]$ will not add any sequences of ordinals of length $< \ka_{n+1}$. 
		
		By Lemma \ref{branchapprox}, it follows that $\dot{b}$ is thinly $\ka_{n+2}$-approximated by the full forcing $\Add(\ka_{n+1}, \pi(\ka_{n+2}))^V \times \PP_1 \times \PP_2 \times \PP_3 \times \Add(\ka_n, \pi(\ka_{n+2}))^V$ over $V[A][U][S][e]$, and thus $\ka_{n+1}$-approximated by the same forcing over $V[A][U][S][e]$. Since the poset $\Add(\ka_{n}, \pi(\ka_{n+3}))^V$ is $\ka_{n+1}$-cc, we can apply Lemma \ref{lem:cc+closedtoapprox} to conclude that $\dot{b}$ is $\ka_{n+1}$-approximated over $V[A][U][S][e][\hat{A}_{n+1}][G_1\times G_2\times G_3]$.
		Thus the branch could not have been added by $\hat{A}_n$. We conclude that $b\in V[A][U][S][e][\hat{A}_{n+1}][G_1\times G_2\times G_3]$.
	\end{proof}
	
	\begin{claim}
		The branch $b$ is in $V[A][U][S][e][\hat{A}_{n+1}]$.
	\end{claim}
	\begin{proof}
		To show that adding $G_1\times G_2\times G_3$ does not add a branch, we follow the argument of \cite[Lemma 4.29]{NeemanTPNw+1}. Recall that $\PP_1\times\PP_2 \times\PP_3$ is $\ka_{n+1}$-closed in $W\defeq V[A][U][S\uhr[\ka_{n+1},\nu)][\hat{A}_{n+1}]$,
		and that $V[A][U][S][e][\hat{A}_{n+1}]$ is an extension of $W$ by a $\ka_{n+1}$-cc poset. Note also that $2^{\ka_n} > \ka_{n+1}$ in $W$. By Lemma \ref{lem:approxcc+closed}, $V[A][U][S][e][\hat{A}_{n+1}][G_1\times G_2\times G_3]$ has the thin $\ka_{n+2}$-approximation property over $V[A][U][S][e][\hat{A}_{n+1}].$ 
		
		Let $\dot{b}$ be a $\PP_1\times\PP_2\times\PP_3$-name for $b$. Note that by Lemma \ref{branchapprox}, $\dot{b}$ is thinly $\ka_{n+2}$-approximated by the forcing $\Add(\ka_{n+1}, \pi(\ka_{n+3}))^V\times (\PP_1\times \PP_2\times \PP_3)$ over $V[A][U][S][e].$ Since $\Add(\ka_{n+1},\pi(\ka_{n+3}))^V$ is $\ka_{n+2}$-cc, we can apply Lemma \ref{lem:cctoapprox} to see that $\dot{b}$ is thinly $\ka_{n+2}$-approximated by $\PP_1\times\PP_2\times\PP_3$ over $V[A][U][S][e][\hat{A}_{n+1}]$.
		It follows that $b$ cannot have been added by $G_1\times G_2\times G_3$, and so $b \in V[A][U][S][e][\hat{A}_{n+1}]$.
	\end{proof}
	
	\begin{claim}
		The branch $b$ is in $V[A][U][S][e]$.
	\end{claim}
	\begin{proof}
		Since the square of $\Add(\ka_{n+1},\pi(\ka_{n+3}))^V$ is $\ka_{n+2}$-cc over $V[A][U][S][e]$, by Lemma \ref{lem:ccbranch} we see that $\Add(\ka_{n+1},\pi(\ka_{n+3}))^V$ has the $\ka_{n+2}$-approximation property over $V[A][U][S][e]$. Let $\dot{b}$ be a $\Add(\ka_{n+1},\pi(\ka_{n+3}))^V$-name for $b$. By Lemma \ref{branchapprox}, $\dot{b}$ is $\ka_{n+2}$-approximated by $\Add(\ka_{n+1},\pi(\ka_{n+3}))^V$ over the model $V[A][U][S][e]$. Thus $\hat{A}_{n+1}$ cannot have added the branch, and so $b \in V[A][U][S][e]$.
	\end{proof}
	
	We have shown that the branch $b$ found in the larger model was actually present in the target model $V[A][U][S][e]$; as described in Lemma \ref{lem:embeddingtobranch}, $b$ must be ineffable.
\end{proof}

\subsection{The strong tree property up to \texorpdfstring{$\aleph_{\w+1}$}{ℵ\_ω+1}}\label{ss:TPatNw+1}

Finally, we put all the pieces together. The argument here is very similar to \cite[Section 6]{NeemanTPNw+1}. We replace \cite[Lemma 3.10]{NeemanTPNw+1} with Theorem \ref{thm:TP2cardfullNA}, and modify some details to accommodate the auxiliary collapse.

\begin{theorem}
	Let $\langle \ka_n\mid 2 \leq n<\w\rangle$ be an increasing sequence of supercompact cardinals with supremum $\nu$. Then there are generics $A$ for $\A$, $U$ for $\UU$ over $V[A]$, $S$ for $\CC^{+A\ast U}$ over $V[A][U]$, and $e$ for $\Coll(\w, \mu)$ over $V[A][U][S]$, such that in the final model $V[A][U][S][e]$ we have the following properties:
	\begin{itemize}
		\item $\aleph_n = \ka_n$ for all $n < \w$
		\item $\aleph_\w$ is strong limit
		\item $(\nu^+)^V = \aleph_{\w+1}$
		\item the strong tree property holds at $\aleph_{\w+1}$
		\item ITP holds at $\aleph_n$ for $2 \leq n < \w$.
	\end{itemize}
\end{theorem}
\begin{proof}
	In order to apply Theorem \ref{thm:TP2cardfullNA}, we need to separate the pieces of the forcing that depend directly on $\ka_1$. In particular, we wish to examine the forcing with $\Add(\ka_1^+, \ka_3)^V \times \CC_0(\ka_1^+)^{+A_0\ast U_0}\times \Coll(\w,\ka_1)$ removed. Unfortunately $\UU_{[1,\w)}$ and $\CC^{+A\ast U}_{[1,\w)}$ rely on the generics for the forcings we wish to remove, so we pass to the versions of these posets that do not depend on the earlier generics. To do this, we need the following definition.
	
	\begin{defn}
		Let $F$ be a filter on $\A\uhr\beta \ast \dot{\UU}\uhr\beta$. Let $\gamma\leq \nu$ and let $B$ be a filter on $\BB^{+F}\uhr[\beta,\gamma)$. We define a filter $F+B$ on $\A\uhr\gamma \ast \dot{\UU}\uhr\gamma$ by $F+B = \{(a,u) \mid (a,u
		\uhr\beta)\in F \text{ and } u\uhr[\beta,\gamma)\in B\}$.
	\end{defn}
	
	Let $A_0\ast U_0$ be generic for $\A_0\ast\dot{\UU}_0$ over $V$. Let $B_{[1,\w)}$ be generic for $\BB_{[1,\w)}^{+A_0\ast U_0}\uhr[\ka_2,\nu)$ over $V[A_0\ast U_0]$. Let $C_{[1,\w)}$ be generic for $\CC^{+A_0\ast U_0+B}\uhr[\ka_2,\nu)$ over $V[A_0\ast U_0][B_{[1,\w)}]$. Let $A_{[2,\w)}$ be generic for $\A_{[2,\w)}$ over $V[A_0\ast U_0][B_{[1,\w)}][C_{[1,\w)}].$ Let $M$ denote the model $V[A_{[2,\w)}][A_0\ast U_0][B_{[1,\w)}][C_{[1,\w)}]$.
	
	For each $\mu \in \Index$, let $\LL_\mu$ be the poset $\Add(\mu^+, \ka_3)^V \times \CC_0(\mu^+)^{+A_0\ast U_0}\times \Coll(\w,\mu)$. $\LL_\mu$ is the remaining piece of the forcing that depends on $\ka_1$, with $\mu^+$ chosen to be the value of $\ka_1$. Recall that $\ka_1$ is a parameter in the definition of $\CC_0$, so by $\CC_0(\mu^+)$ we mean $\CC_0$ defined relative to the parameter $\ka_1 = \mu^+$.
	
	Let $\la \geq \nu^+$ such that $\la^\nu =\la$. Let $K$ be generic for $\Coll(\nu^+, \la)^{V}$. Since $K$ is mutually generic with the remaining pieces of the forcing, we have that $M[K] = V[K][A_{[2,\w)}][A_0\ast U_0][B_{[1,\w)}][C_{[1,\w)}]$.
	
	\begin{claim} \label{claim:tpcondition1}
	In $V[K][A_{[2,\w)}]$, there is a $\nu^+$-supercompactness embedding $\pi:V[K][A_{[2,\w)}] \to V^*[K^*][A^*_{[2,\w)}]$ with critical point $\ka_2$ and $|\pi(\ka_2)| = \nu^{++}$, such that $\nu \in \pi(\Index)$. In any extension $M[K][\hat{A}_0]$ of $M[K]$ by the poset $\Add(\w, [\ka_2, \pi(\ka_2)))^V$, $\pi$ extends to an elementary embedding $\pi:M[K] \to M^*[K^*]$ with $\nu \in \pi(\Index)$.
	\end{claim}
	\begin{proof}
	The existence of this embedding in $V[K][A_{[2,\w)}]$ is immediate from the indestructibility of $\ka_2$, noting that $K\times A_{[2,\w)}$ is generic for a $\ka_2$-directed closed forcing. The fact that this embedding extends to have domain $M[K]$ and the fact that $\nu \in \pi(\Index)$, are analogous to the proof of \cite[Lemma 5.7]{NeemanTPNw+1}.
	To verify that $\nu \in \pi(\Index)$, we can apply the arguments of \cite[Lemma 5.7]{NeemanTPNw+1} in $M[K] = V[K][A_{[2,\w)}][A_0\ast U_0][B_{[1,\w)}][C_{[1,\w)}]$; since $\Coll(\nu^+, \la)^V$ is $\nu^+$-closed in $V$, there is no obstacle to doing so.
	\end{proof}
	
	The following claim is a strengthening of \cite[Lemma 5.8]{NeemanTPNw+1}, and is why we require the auxiliary collapse $\Coll(\nu^+, \la)^{V}$. In order to have the small models contain all relevant objects, they need to be of size at least $|\la|$; since we have collapsed $\la$ to $\nu^+$, however, the $\nu^+$-closure of our posets is sufficient to obtain the desired generic.
	\begin{claim}\label{claim:tpcondition1part2}
		 Let $R$ be a rank initial segment of the universe, large enough to contain all relevant objects. Let $\bar{M}[K] = \bar{V}[K][A_{[2,\w)}][A_0\ast U_0][B_{[1,\w)}][C_{[1,\w)}]$, where $\bar{V}$ is the transitive collapse of $X \prec R$ with $X \in V, V_\nu \subseteq X, |X| = \la$, $\la \subseteq X$
		 , and $X$ is closed under sequences of length $\nu$ in $V$. Let $\hat{A}_0$ be generic for $\Add(\w, (\nu^{++})^{\bar{M}})$ over $M[K]$ (and also $\bar{M}[K]$). Let $\pi:\bar{M}[K] \to \bar{M}^*[K^*]$ be the embedding from the previous claim applied in $\bar{M}[K][\hat{A}_0]$. Let $e$ be generic for $\Coll(\w, \nu)$ over $M[K][\hat{A}_0]$.
		 
		Then in $M[K][\hat{A}_0][e]$ there are filters $A^*_1$ and $S^*_0$ so that $A^*_1\times S^*_0\times e$ is generic for $\pi(\LL)(\nu)$ over $\bar{M}^*[K^*]$.
	\end{claim}
\begin{proof}
	As in \cite[Lemma 5.8]{NeemanTPNw+1}, it is enough to find a generic $A_1^*\times C_0^*$ for the poset $\Add(\nu^+, \pi(\ka_3))^{\bar{V}^*} \times \pi(\CC_0)(\nu^+)$ over $\bar{V}[K][A_{[2,\w)}]$ that belongs to $V[K][A_{[2,\w)}]$. This can be done since $\Add(\nu^+, \pi(\ka_3))^{\bar{V}^*} \times \pi(\CC_0)(\nu^+)$ is $\nu^+$-closed in $V[K][A_{[2,\w)}]$, and $\bar{V}[K][A_{2,\w)}]$ has size $\nu^+$.
\end{proof}
	Note that since $\la \subseteq X$, $\Coll(\nu^+, \la)^V$ is in the transitive part of $X$, and thus is sent to itself by the transitive collapse.
	\begin{claim}\label{claim:tpcondition2}
		In $M[K]$, for all $m \geq 3$, there is a generic $\nu^+$-supercompactness embedding with critical point $\ka_m$. The poset adding this embedding is the product of $\Add(\ka_m, \pi(\ka_{m+2}))^V$, $\Add(\ka_{m+1}, \pi(\ka_{m+3}))^V$, and $\pi(\CC)^{+F\uhr\ka_{m+2}}\uhr(\ka_{m+2}, \pi(\ka_{m+2}))$, where $F$ is defined to be $A_0\ast U_0^{+B_{[1,\w)}}.$ The full support $\ka_{m}$-th power of this poset is $<\ka_m$-distributive over $M[K]$.
	\end{claim}
	\begin{proof}
		The proof is the argument of \cite[Lemma 5.6]{NeemanTPNw+1}, carried out over $V[K]$ rather than $V$. Since $K$ is generic for a $\nu^+$-directed-closed poset, each $\ka_n$ remains indestructibly supercompact in $V[K]$, and all necessary properties of the other posets over $V$ will still hold in $V[K]$.
	\end{proof}

	\begin{claim}\label{cl:TPnu+}
		There is $\mu \in \Index$ so that in the extension of $M$ by $\LL_\mu$, the strong tree property holds at $\nu^+$.
	\end{claim}
	\begin{proof}
		It suffices to check that the hypotheses of Theorem \ref{thm:TP2cardfullNA} hold for all $\la$, with $M$ serving as our ground model. First we note that $M\models \ka_n^+ = \ka_{n+1}$ for all $n\geq 2$, and $\LL_\mu$ is $\ka_2$-cc.
		The embedding for condition (1) is obtained from Claims \ref{claim:tpcondition1} and \ref{claim:tpcondition1part2}, noting that for each $\la$, the poset $\Add(\w, [\ka_2, \pi(\ka_2)))^V$ adding the embedding is the product of a $\nu^+$-Knaster poset and a $\nu^+$-closed poset. The embeddings for condition (2) are obtained from Claim \ref{claim:tpcondition2}.
	\end{proof}
	
	Let $\mu$ be given by the previous claim, and let $A_1 \times S_0 \times e$ be generic for $\LL_\mu$ over $M$. Let $A = A_0\times A_1\times A_{[2,\w)}$. Let $U_{[1,\w)}$ be the upwards closure of $B_{[1,\w)}$ in $\UU_{[1,\w)}$, and let $U = U_0 \times U_{[1,\w)}$. Let $S_{[1,\w)}$ be the upwards closure of $C_{[1,\w)}$ in $\CC^{+A\ast U}_{[1,\w)}$ and let $S = S_0 \ast S_{[1,\w)}$. Let $N$ be the model $V[A][U][S]$.
	
	\begin{prop}\cite[Lemma 5.9]{NeemanTPNw+1}\label{prop:TPNeemanLem5.9}
		There is a $\mu^+$-closed extension $N[G]$ containing $M[A_1\times S_0]$, with $G$ still generic over $N[e]$ and both $\nu$ and $\nu^+$ still cardinals in $N[e][G]$.
	\end{prop}
	\begin{proof}
	The relevant model is the extension of $N$ by the product of the factor forcing refining $U_{[1,\w)}$ to a filter for $\BB^{+A_0\ast\dot{\UU}_0}\uhr[\ka_2,\nu)$ and the factor forcing refining $S_{[1,\w)}$ to a filter for $\CC^{+A_0\ast\dot{\UU}_0}\uhr[\ka_2,\nu)$. By Fact \ref{fact:Brefinementisdirclosed} and Fact \ref{fact:Crefinementdirclosed}, every descending sequence of conditions in these posets with length $<\ka_2$ belonging to $V[A_0\ast U_0]$ have lower bounds. From Fact \ref{fact:laterforcingsdist} we see that every descending sequence of conditions in these posets with length less than $\ka_1 = \mu$ in $V[A][U][S]$ belong to $V[A\uhr\ka_2][U\uhr\ka_1]$, so they must belong to $V[A_0\ast U_0]$. We conclude that the factor posets are $\mu$-closed in $V[A][U][S]$. When we extend $N[e]$ by these posets, the resulting model is contained in $V[A][U\uhr\ka_1][B\uhr[\ka_1,\nu)][C][e]$; this model preserves $\nu$, $\nu^+$, and each $\ka_n$. (For more details about the cardinal preservation, see \cite[Lemma 4.24]{NeemanTPNw+1} and \cite[Remark 4.25]{NeemanTPNw+1}.)
	\end{proof}

	In particular, $N[G]$ projects to $M[A_1\times S_0]$, which projects to $N$. Moreover $M[\LL_\mu] = M[A_1\times S_0 \times e]$, which projects to our final model $N[e]$.
	
	\begin{claim}\label{cl:TPnu^+inN[e]}
		In $N[e]$, the strong tree property holds at $\nu^+$.
	\end{claim}
	\begin{proof}
		This is almost identical to the proof of \cite[Claim 6.5]{NeemanTPNw+1}, using a slightly more general branch lemma.
		Let $d$ be a thin $\mc{P}_{\nu^+}(\la)$-list. By Claim \ref{cl:TPnu+}, $d$ has a cofinal branch in the model $M[A_1\times S_0 \times e]$.
		
		Let $N[G]$ be the extension of $N$ given by Proposition \ref{prop:TPNeemanLem5.9}.
		Since $\dot{d}$ has a cofinal branch in $M[A_1\times S_0\times e]$, this branch is also present in $N[e][G]$. Applying Lemma \ref{lem:M-S} and Lemma \ref{branchapprox}, we conclude that the branch could not have been added by $G$. It follows that $\dot{d}$ has a cofinal branch in $N[e]$ as desired.
	\end{proof}
	
	\begin{claim}
		In $N[e]$, ITP holds at $\aleph_{n}$ for all $n\geq 2$.
	\end{claim}
	\begin{proof}
		This follows immediately from Theorem \ref{thm:ITPatNn}, noting that $\ka_n = \aleph_n$ in $N[e]$.
	\end{proof}
	
	We conclude that in $N[e]$, ITP holds at $\aleph_n$ for $2 \leq n < w$, and the strong tree property holds at $\nu^+ = \aleph_{\w+1}$. This completes the proof.
\end{proof}

\section{ITP at the successor of a singular cardinal}\label{s:ITPsuccsing}
In this section, we give an analogue of Theorem \ref{thm:TP2cardfullNA}, describing a general class of forcings which will obtain ITP at the successor of a singular cardinal. In particular, we show that ITP can be obtained at the successor of a singular cardinal with uncountable cofinality. Note that the hypotheses of this theorem are a little bit stricter; we impose more constraints on the structure of $\LL_\mu$, and we require $\ka_0$ to be supercompact rather than merely generically supercompact. Our argument is a generalization of the techniques in \cite{HachtmanITPNw+1}.

As is becoming standard, we will prove the one-cardinal case separately for clarity.

\subsection{The One-Cardinal Case}

\begin{theorem}\label{thm:gen1cardITP}
	Let $\tau$ be a regular cardinal. Let $\langle \ka_\rho \mid \rho < \tau \rangle$ be an increasing continuous sequence of cardinals above $\tau$ with supremum $\nu$, such that $\ka_{\rho}^+ = \ka_{\rho+1}$ for all $\rho < \tau$. Let $\Index$ be a subset of $\ka_0$, and fix $\rho' < \tau$. For each $\mu \in \Index$, let $\LL_\mu$ be the product of forcings $\PP_\mu$ and $\Q_\mu$ where $|\PP_\mu| < \mu^+$ and $\Q_\mu$ is $\mu^{++}$-closed, such that $|\LL_\mu| \leq \ka_{\rho'}$.
	In addition, suppose that we have the following:
	\begin{itemize}
		\item $\ka_0$ is $\nu^+$-supercompact, with a normal measure $U_0$ on $\mc{P}_{\ka_0}(\nu^+)$ and corresponding embedding $i$, such that $\nu \in i(\Index)$.
		\item For all ordinals $\rho < \tau$ there is a generic $\nu^+$-supercompactness embedding $j_{\rho+2}$ with domain $V$ and critical point $\ka_{\rho+2}$, added by a poset $\FF$ such that the full support power $\PP^{\ka_\rho}$ is $<\ka_{\rho}$-distributive in $V$.
	\end{itemize}
	Then there exists $\mu \in \Index$ such that $\ITP(\nu^+, \nu^+)$ holds in the extension of $V$ by $\LL_\mu$.
\end{theorem}
\begin{proof}
	Suppose not. For each $\mu \in \Index$, let $\dot{d}^\mu$ be a $\LL_\mu$-name for a thin $\mc{P}_{\nu^+}(\nu^+)$-list with no ineffable branch. Assume that the $\alpha$-th level of $\dot{d}^\mu$ is enumerated by the names $\{\dot{\sigma}_\alpha^\mu(\xi) \mid \xi < \nu\}$, and that (for sufficiently large $\alpha$) there are no repetitions in this sequence.
	
	By assumption, we have a normal measure $U_0$ on $\mc{P}_{\ka_0}(\nu^+)$, with corresponding embedding $i : V \to M$. Let $\ka_x$ denote $\sup(\ka_0 \cap x)$. Recall that $\ka_x = \ka_0 \cap x$ on a measure one set, and $[x \mapsto \ka_x]_{U_0} = \ka_0$. Then $[x \mapsto \ka_x^{+\tau}]_{U_0} = \ka_0^{+\tau} = \nu$. Let $\mu_x$ denote $\ka_x^{+\tau}$. Note that since $\nu \in i(\Index)$, $\mu_x \in \Index$ on a measure one set. Note also that $i\dot{d}_{\sup i''\nu^+}^{\nu} = [x \mapsto \dot{d}_{\sup x}^{\mu_x}]_{U_0}$.

\begin{lemma}\label{lem:ITP1cardpigeonhole}
	There exists a successor ordinal $\rho$ with $\rho' < \rho < \tau$, an unbounded $S \subseteq \nu^+$, $A \in U_0$, and a map $x \mapsto (p_x, q_x)$ such that for all $x \in A$ and $\alpha \in x\cap S$, there is $\xi < \ka_\rho$ such that $(p_x, q_x) \forces_{\LL_{\mu_x}} \dot{d}_{\sup x}^{\mu_x} \cap \alpha = \dot{\sigma}^{\mu_x}_\alpha(\xi)$.
\end{lemma}
\begin{proof}
		Let $\LL$ denote $[x \mapsto \LL_{\mu_x}]_{U_0}$. By assumption, each $\LL_{\mu_x}$ is the product $\PP_{\mu_x} \times \Q_{\mu_x}$, where $|\PP_{\mu_x}| < \mu_x^+$ and $\Q_{\mu_x}$ is $\mu_x^{++}$-closed. Let $\PP = [x \mapsto \PP_{\mu_x}]_{U_0}$ and $\Q = [x\mapsto \Q_{\mu_x}]_{U_0}$. We conclude that $|\PP| < \nu^+$, $\Q$ is $\nu^{++}$-closed, and $\LL = \PP \times \Q$.
		
		For all $\alpha < \nu^+$, there is some successor ordinal $\rho_\alpha < \tau$, $\xi < i(\ka_{\rho_\alpha})$, and $(p_\alpha, q_\alpha) \in \LL$ such that $(p_\alpha, q_\alpha) \forces_{\LL} i\dot{d}^{\nu}_{\sup i''\nu^+}\cap i(\alpha) = i\dot{\sigma}^{\nu}_{i(\alpha)}(\xi).$ Since $\Q$ is $\nu^{++}$-closed, working inductively we can choose the conditions $q_\alpha$ to be decreasing, with lower bound $q$. Since $|\PP| < \nu^+$ and $\tau < \nu^+$, there is an unbounded $S \subseteq \nu^+$, a fixed successor ordinal $\rho < \tau$, and a fixed $p \in \PP$ such that for all $\alpha \in S$, $\rho_\alpha = \rho$ and $p_\alpha = p$. Let $[x \mapsto p_x]_{U_0} = p$ and $[x \mapsto q_x]_{U_0} = q$.
		
		Applying \L os' theorem, we conclude that for all $\alpha \in S$, there is a measure one set $A_\alpha$ such that for all $x \in A_\alpha$, there is $\xi < \ka_\rho$ such that $(p_x, q_x) \forces_{\LL_{\mu_x}} \dot{d}^{\mu_x}_{\sup x} \cap \alpha = \dot{\sigma}_\alpha^{\mu_x}(\xi)$. Let $A \defeq \triangle_{\alpha \in S} A_\alpha$. This is a measure one set with the desired properties.
\end{proof}	

Let $I = \{(\mu,p,q) \mid \mu \in \Index, (p,q) \in \LL_\mu\}$. For all $s = (\mu,p,q)$ in $I$, we define the relation $R_s$ on $S \times \ka_\rho$ by $(\alpha, \zeta) R_s (\beta, \xi)$ iff $\alpha \leq \beta$ and $(p,q) \forces \dot{\sigma}^\mu_\alpha (\zeta) = \dot{\sigma}_\beta^\mu(\xi)\cap \alpha$.

\begin{lemma}
	$\langle R_s \rangle_{s\in I}$ is a system on $S \times \ka_\rho$.
\end{lemma}
\begin{proof}
	The first two conditions are trivial. For the third, let $\alpha < \beta$ both in $S$, and let $x \in A$ such that $\alpha, \beta \in x$. Then there exists  $(p_x, q_x) \in \LL_{\mu_x}$ and $\zeta, \xi < \ka_\rho$ such that $(p_x, q_x) \forces \dot{\sigma}_\alpha^{\mu_x}(\zeta) = \dot{\sigma}_\beta^{\mu_x}(\xi)\cap \beta$. In particular, letting $s = (\mu_x, p_x, q_x)$, we see that $(\alpha, \zeta) R_s (\beta, \xi)$ as desired.
\end{proof}

\begin{lemma}\label{lem:ITP1cardgettingabranch}
	There exists an unbounded $S' \subseteq S$ and a system of branches $\langle b_{s,\delta} \mid s\in I, \delta < \ka_\rho \rangle$ through $\langle R_s \uhr S'\times \ka_\rho\rangle_{s\in I}$ such that each $b_{s,\delta}$ is a branch through $R_s\uhr S'\times \ka_\rho$.
\end{lemma}
\begin{proof}
By assumption, we have a generic $\nu^+$-supercompactness embedding $j$ with critical point $\ka_{\rho+3}$, added by a poset $\FF$ with generic $F$ whose full support power $\FF^{\ka_{\rho+1}}$ is $<\ka_{\rho+1}$-distributive in $V$. Thus $V$ satisfies hypothesis (2) of Lemma \ref{lem:1cardbranch}. Note that $\FF$ itself must also be $<\ka_{\rho+1}$-distributive in $V$. Work in $V[F]$.

Let $\gamma \in j(S) \sm \sup j''\nu^+$. Note that each $\LL_\mu$ has size $<\ka_{\rho+3}$. As before, we can assume that the underlying set of $\LL_\mu$ is an ordinal, so that $\LL_\mu$ is a bounded subset of $V_{\ka_{\rho+3}}$; it follows that $j(I) = I$. For each $\delta < \ka_\rho$ and $s = (\mu, p,q) \in I,$ we define
\[b_{s,\delta} = \{(\alpha, \zeta) \mid \alpha \in S, \zeta < \ka_\rho, (p,q) \forces_{\LL_\mu} j(\dot{\sigma}_\alpha^\mu(\zeta)) = j\dot{\sigma}_\gamma^\mu(\delta) \cap j(\alpha).\}\]

We claim that $\langle b_{s,\delta} \mid s\in I, \delta < \ka_\rho\rangle$ is a system of branches through $\langle R_s\rangle_{s\in I}$. By construction, each $b_{s,\delta}$ is linearly ordered and downwards closed. It remains to verify that $\bigcup \dom(b_{s,\delta}) = S$.

Since $\ka_\rho < \crit(j)$, we can apply elementarity to Lemma \ref{lem:ITP1cardpigeonhole}, concluding that for all $\alpha \in S$ and $x \in j(A)$ with $j(\alpha), \gamma$ both in $x$, there exists $\zeta, \delta < \ka_\rho$ and $s = (\mu_x, p_x, q_x) \in j(I) = I$ such that
\[(p_x,q_x) \forces_{\LL_{\mu_x}} j(\dot{\sigma}_\alpha^{\mu_x}(\zeta)) = j\dot{\sigma}_\gamma^{\mu_x}(\delta) \cap j(\alpha).\]
In particular, $\bigcup \dom(b_{s,\delta}) = S$. This system may not belong to $V$, since it is defined in $V[F]$, but it satisfies condition (1) of Lemma \ref{lem:1cardbranch}. Applying that lemma, we conclude that there is some $(s, \delta) \in I \times \ka_{\rho}$ such that $b_{s,\delta}$ is cofinal and belongs to $V$.

Let $\mc{D} = \{(s, \delta)\mid b_{s,\delta}\in V\}$. Since $\FF$ is $<\ka_{\rho+1}$-distributive, $\mc{D} \in V$, and in fact $\langle b_{s,\delta} \mid (s,\delta) \in \mc{D}\rangle$ is also in $V$. Since $\mc{D}$ must contain at least one pair $(s,\delta)$ corresponding to a cofinal branch, the set $S' = \bigcup_{(s,\delta) \in \mc{D}} \dom(b_{s,\delta})$ is unbounded in $\nu^+$, and $\langle b_{s,\delta}\rangle_{(s,\delta)\in \mc{D}}$ is a system of branches through $\langle R_s \uhr S' \times \ka_\rho\rangle_{s\in I}$.

In addition, by taking a subset of $I \times \ka_\rho$ if needed, we may assume that for all $s \in I$ and $\eta < \delta < \ka_\rho$, if $b_{s,\eta}$ and $b_{s,\delta}$ are both cofinal, then they are distinct. In particular, we can assume that if $b_{s,\delta}$ and $b_{s,\eta}$ are (defined and) equal cofinally often, then $\delta = \eta$. We can do this by removing any duplicates, which may appear if $j\dot{\sigma}_\gamma^\mu(\eta)$ and $j\dot\sigma_\gamma^\mu(\delta)$ are above $\sup j''\nu^+$. (Note that by the distributivity of $\FF$, the necessary subset of $I \times \ka_\rho$ will still be in $V$.)
\end{proof}

For all $(s, \delta) \in \mc{D}$, we define $\dot{\pi}_{s,\delta} = \bigcup\{\dot{\sigma}_\alpha^\mu(\zeta) \mid (\alpha, \zeta) \in b_{s,\delta}\}$. As in \cite{HachtmanITPNw+1}, these branches have some useful properties. Let $s = (\mu,p,q)$ and $s' = (\mu, p', q')$. If $(p',q')\leq (p,q)$, then $R_s \subseteq R_{s'}$, and $b_{s,\delta} \subseteq b_{s',\delta}$ for all $\delta < \ka_\rho$. Moreover, if $b_{s,\delta}$ is cofinal and $(s,\delta)\in \mc{D}$, then $(s', \delta) \in \mc{D}$. If $b_{s,\delta}$ is cofinal, then $(p,q)$ forces that $\dot{\pi}_{s,\delta}$ is a cofinal branch through $\dot{d}^\mu$.

Next, we wish to bound the splitting for all branches (including those not in $V$). Working in $V[F]$, for each $\eta < \delta < \ka_\rho$ and $s\in I$, we define $\alpha_{s,\eta,\delta}$ as follows. If $b_{s,\delta}$ and $b_{s,\eta}$ are both bounded, we define $\alpha_{s,\eta,\delta}$ to be $\sup \left(\dom(b_{s,\eta}) \cup \dom(b_{s,\delta})\right)$. If not, then we define $\alpha_{s,\eta,\delta}$ to be the least $\alpha$ such that for all $\alpha' > \alpha$, $b_{s,\eta}(\alpha')$ and $b_{s,\delta}(\alpha')$ are not both defined and equal. Note that each $\alpha_{s,\eta,\delta}$ is below $\nu^+$; if it isn't, then $b_{s,\eta,\delta}$ must agree on cofinally many $\alpha$, but since we have removed duplicate branches this would mean that $\delta = \eta$.

Let $\bar{\alpha} = \sup_{s\in I,\eta<\delta<\ka_\rho} \alpha_{s,\eta,\delta}+1$. Then if $\alpha > \bar{\alpha}$ and $b_{s,\delta}(\alpha) = b_{s,\eta}(\alpha)$, $\delta = \eta$ and $b_{s,\delta}$ is cofinal.

For all $x \in A$, let $(p_x, q_x) \in \LL_{\mu_x}$ be as in Lemma \ref{lem:ITP1cardpigeonhole}. Define $s_x = (\mu_x, p_x, q_x)$.

\begin{lemma}\label{lem:1cardITPdagger}
	There exists an unbounded $\bar{S}\subseteq S'$ and $\bar{A}\in U_0$ with $\bar{A} \subseteq A$ such that for all $x \in \bar{A}$ and all $\alpha \in \bar{S}\cap x$, the following statement holds:
	\[(\dagger_{x,\alpha}) \quad \exists \delta < \ka_n \quad (s_x,\delta) \in \mc{D} \text{ and } (p_x,q_x) \forces_{\LL_{\mu_x}} \dot{d}^{\mu_x}_{\sup x} \cap \alpha = \dot{\pi}_{s_x,\delta}\cap \alpha. \]
\end{lemma}
\begin{proof}
	Define $A_\alpha \defeq \{x\in A \mid (\dagger_{x,\alpha}) \text{ holds }\}$. It suffices to show that $\bar{S} = \{\alpha \in S' \mid A_\alpha \in U_0\}$ is unbounded in $\mu$; if this holds, then $\bar{A} = A \cap \triangle_{\alpha \in \bar{S}} A_\alpha$ will have the desired properties.
	
	Suppose $\bar{S}$ is bounded. Fix $\alpha_0 < \nu^+$ such that $\bar{\alpha} < \alpha_0$ and $A_\alpha \notin U_0$ for all $\alpha > \alpha_0$ that are in $S'$. Define $A' \defeq A \cap \triangle_{\alpha_0 < \alpha \in S'} \mc{P}_{\ka_0}(\nu^+)\sm A_\alpha$. Then $A' \in U_0$, and $\dagger_{x,\alpha}$ fails whenever $\alpha_0 < \alpha$, $\alpha \in x\cap S'$, and $x \in A'$.
	
	Let $R'_s$ be obtained by removing every ground model branch from $R_s$. That is, for each $s\in I$, $(\alpha, \zeta) R'_s (\beta, \xi)$ if and only if $\alpha_0 < \alpha,\beta \in S', (\alpha, \zeta) R_s (\beta, \xi)$, and for all $\delta < \ka_\rho$, if $(s,\delta) \in \mc{D}$, then $(\alpha, \zeta) \notin b_{s,\delta}$. We wish to show that $\langle R'_s\rangle_{s \in I}$ is a system on $(S' \sm \alpha_0)\times \ka_\rho$.
	
	As before, the first two conditions are immediate from the definition. For the third, suppose $\alpha_0<\alpha<\beta$ with $\alpha, \beta$ both in $S'$. We need to show that there exists $\zeta, \xi < \ka_{\rho}$ and $s \in I$ such that $(\alpha, \zeta) R'_s (\beta, \xi)$. Let $\gamma \in j(S) \sm \sup j''\nu^+$ be the element used in the previous lemma to define the system of branches. Note that $\{x \mid j(\alpha), j(\beta), \gamma \in x\}$ is a club in $\mc{P}_{\ka_0}(j(\nu^+))$. Since $A'$ intersects every club in $\mc{P}_{\ka_0}(\nu^+)$, by elementarity we conclude that there exists some $x' \in j(A')$ with $j(\alpha), j(\beta)$, and $\gamma$ all in $x'$. Note that $j(\alpha), j(\beta) \in j(S)$. Applying elementarity to Lemma \ref{lem:ITP1cardpigeonhole}, we conclude that there exist $\zeta, \xi, \delta < \ka_\rho$ such that
	\[(p_{x'},q_{x'}) \forces j\dot{\sigma}_{j(\alpha)}^{\mu_{x'}}(\zeta) = j\dot{\sigma}_{\gamma}^{\mu_{x'}}(\delta)\cap j(\alpha), j\dot{\sigma}_{j(\beta)}^{\mu_{x'}}(\xi) = j\dot{\sigma}_\gamma^{\mu_{x'}}(\delta)\cap j(\beta),\]
	and each of these are forced by $(p_{x'}, q_{x'})$ to cohere with $\dot{j}d_{\sup x}^{\mu_x}$. 
	Note that $x'$ may not be the image of an element from the ground model. Since $s_{x'} \in j(I) = I$, however, by elementarity there exists $x\in A'$ with $\alpha, \beta \in x$ such that $s_x = s_{x'}$. We conclude that $(\alpha, \zeta)$ and $(\beta, \xi)$ are both in $b_{s_{x'},\delta}$, with $(\alpha, \zeta) R_{s_{x'}} (\beta, \xi)$.
	
	Since we are above the splitting, if $(\alpha, \zeta) \in b_{s_x,\delta'}$ for any $\delta' \neq \delta$, this branch must coincide with $b_{s_x, \delta}$. So to finish showing that $(\alpha, \zeta) R'_{s_x} (\beta, \xi)$, we simply need to show that $(s_x, \delta) \notin \mc{D}$. We will do this by contradiction.
	
	Suppose $(s_x, \delta) \in \mc{D}$. Since $\alpha_0 <\alpha$ and $\alpha \in x\cap S'$ with $x \in A'$, $(\dagger_{x,\alpha})$ must fail. Then
	\[(p_x, q_x) \not\forces \dot{d}_{\sup x}^{\mu_x} \cap \alpha = \dot{\pi}_{s_x,\delta}\cap \alpha. \]
	We have chosen $x$ so that
	\[(p_x, q_x) \forces \dot{d}_{\sup x}^{\mu_x}\cap \alpha = \dot{\sigma}^{\mu_x}_\alpha(\zeta).\]
	But since $(\alpha, \zeta) \in b_{s_x,\delta}$, by the definition of $\pi_{s_x,\delta}$ we must have
	\[(p_x,q_x) \forces \dot{\pi}_{s_x,\delta}\cap \alpha = \dot{\sigma}^{\mu_x}_\alpha(\zeta).\]
	This gives a contradiction. It follows that $(s_x, \delta) \notin \mc{D}$, concluding our proof that $\langle R'_s\rangle_{s\in I}$ is a system.
	
	For each $(s,\delta) \notin \mc{D}$, let $b'_{s,\delta}$ be the restriction of $b_{s,\delta}$ to $R'_s$. Then $\langle b'_{s,\delta} \mid (s,\delta) \notin \mc{D}\rangle$ is a system of branches through $\langle R'_s\rangle_{s \in I}$. Repeating the argument of Lemma \ref{lem:ITP1cardgettingabranch}, we conclude that there exists some $s \in I$ and $\delta < \ka_\rho$ such that $b'_{s,\delta} \uhr R_s'$ is cofinal and belongs to $V$. Since we can recover $b_{s,\delta}$ from any cofinal subset, we conclude that $(s,\delta)$ must be in $\mc{D}$, contradicting our definition of $R_s'$. This contradiction proves our initial claim that $\bar{S}$ is unbounded; thus we can construct $\bar{A}$ with the desired properties as described above.
\end{proof}

Let $S^* = \bar{S} \sm (\bar{\alpha}+1)$; this set is still unbounded. Let $A^* = \{x \in A \mid x \cap S^* \text{ cofinal in } \sup x\}$, noting that $A^* \in U_0$. For all $x \in A^*$ and $\alpha \in x\cap S^*$, the witness $\delta$ to $(\dagger_{x,\alpha})$ depends only on $x$: if we have $\alpha < \beta$ both in $S^*$ and $\delta, \delta' < \ka_n$, where $(p_x,q_x) \forces \dot{d}^{\mu_x}_{\sup x}\cap \alpha = \dot{\pi}_{s_x, \delta}\cap \alpha$ and $(p_x,q_x) \forces \dot{d}^{\mu_x}_{\sup x}\cap \beta = \dot{\pi}_{s_x, \delta'}\cap \beta$, then clearly $(p_x,q_x) \forces \dot{\pi}_{s_x, \delta}\cap \alpha = \dot{\pi}_{s_x, \delta'}\cap \alpha$. Since the branches are forced to cohere up to $\alpha$, and $\alpha$ is above the splitting, then $(p_x, q_x)$ must force them to be equal; since we have removed duplicate branches, it follows that $\delta = \delta'$.

Since the witness $\delta$ depends only on $x$, we see that for all $\alpha \in S^*\cap x$, $(p_x, q_x) \forces \dot{d}_{\sup x}^{\mu_x}\cap \alpha = \dot{\pi}_{s_x, \delta}\cap \alpha$. We conclude that $(p_x, q_x) \forces \dot{d}_{\sup x}^{\mu_x} = \dot{\pi}_{s_x, \delta}\cap \sup x.$

For each $s = (\mu,p,q)$ and $\delta$ such that $(s,\delta) \in \mc{D}$, we define $T_{s,\delta} = \{\alpha < \nu^+ \mid (p,q) \forces_{\LL_{\mu}} \dot{d}^{\mu}_\alpha = \dot{\pi}_{s,\delta}\cap \alpha\}$. Let $T = \bigcup_{(s,\delta) \in \mc{D}} T_{s,\delta}$. We have shown that $\{\sup x \mid x \in A^*\} \subseteq T$, so $T$ must be stationary. Since $|\mc{D}| \leq \ka_\rho < \nu^+$, there must be some fixed $(s,\delta)$ such that $T_{s,\delta}$ is stationary. Since $\LL_\mu$ has the $\nu^+$-cc, it will preserve stationary subsets of $\nu^+$. We conclude that $b_{s,\delta}$ defines an ineffable branch through $\dot{d}^\mu$ in any generic extension of $V$ by $\LL_\mu$ containing $(p,q)$.
\end{proof}

\subsection{The Two-Cardinal Case}

\begin{theorem}\label{thm:gen2cardITP}
	Let $\tau$ be a regular cardinal. Let $\langle \ka_\rho \mid \rho < \tau \rangle$ be an increasing continuous sequence of cardinals above $\tau$ with supremum $\nu$, such that $\ka_{\rho}^+ = \ka_{\rho+1}$ for all $\rho < \tau$. Let $\Index$ be a subset of $\ka_0$ such that every $\mu \in \Index$ has cofinality $\tau$. For each $\mu \in \Index$, let $\LL_\mu$ be the product of forcings $\PP_\mu$ and $\Q_\mu$, where $|\PP_\mu| < \mu^+$ and $\Q_\mu$ is $\mu^{++}$-closed, such that $|\LL_\mu| < \ka_{\rho'}$ for some fixed $\rho' < \tau$.
	In addition, suppose that we have the following:
	\begin{itemize}
		\item $\ka_0$ is indestructibly $\nu^+$-supercompact, with a normal measure $U_0$ on $\mc{P}_{\ka_0}(\nu^+)$ and corresponding embedding $i$ such that $\nu \in i(\Index)$.
		\item For all ordinals $\rho < \tau$ and all $\la \geq \nu^+$, there is a generic $\la$-supercompactness embedding $j_{\rho+2}$ with domain $V$ and critical point $\ka_{\rho+2}$, added by a poset $\FF$ such that the full support power $\FF^{\ka_\rho}$ is $<\ka_{\rho}$-distributive in $V$.
	\end{itemize}
	Then there exists $\mu \in \Index$ such that $\ITP$ holds at $\nu^+$ in the extension of $V$ by $\LL_\mu$.
\end{theorem}
\begin{proof}
	As before, let $\ka$ denote $\ka_0$.
	Suppose the theorem fails. Then for every $\mu < \ka$, there is some $\la$ such that $\ITP(\nu^+, \la)$ fails in the extension by $\LL_\mu$. By taking a supremum, we can assume that $\la$ is the same for all $\mu$, and that $\la^{\nu^+} = \la$. For each $\mu \in \Index$, let $\dot{d}^\mu$ be a name for a thin $\mc{P}_{\nu^+}(\la)$ list which is forced by $1_{\LL_\mu}$ not to have an ineffable branch. Let $K$ be generic for $\Coll(\nu^+, \la)^{V}$ over $V$. Note that in $V[K]$, $\la$ is no longer a cardinal, since it has been collapsed; it is simply a set with size $\nu^+$. Since we can assume that $\la$ was regular in $V$, we can likewise assume that $\la$ has cofinality $\nu^+$ in $V[K]$. It is still meaningful to discuss ineffable branches through thin $\mc{P}_{\nu^+}(\la)$ lists when $\la$ is a set rather than a cardinal.
	
	\begin{theorem}\label{thm:genITPinK}
		There is $\mu \in \Index$ such that in $V[K][\LL_\mu]$, $d^\mu$ has an ineffable branch.
	\end{theorem}
	Assuming this theorem, we will complete the proof of Theorem \ref{thm:gen2cardITP}. Let $b$ be an ineffable branch for $d^\mu$ in $V[K][\LL_\mu]$. Since $\Coll(\nu^+, \la)$ is $\nu^+$-closed in $V$, and $\LL_\mu$ is $\ka_\rho$-cc for some $\rho < \tau$, we conclude by Lemma \ref{lem:M-S} that $\Coll(\nu^+, \la)$ is forced to have the thin $\nu^+$-approximation property. Therefore $b \in V[\LL_\mu]$. Since stationarity is downwards absolute, $b$ is an ineffable branch for $d^\mu$ in $V[\LL_\mu]$, contradicting our assumption.
\end{proof}

We will now prove Theorem \ref{thm:genITPinK}. Note that in $V[K]$, $\la$ is no longer a cardinal, having been collapsed to have cardinality and cofinality $\nu^+$.
Since $\la$ has cardinality and cofinality $\nu^+$, $\mc{P}_{\nu^+}(\la)$ is order-isomorphic to $\mc{P}_{\nu^+}(\nu^+)$. Using this isomorphism, we can identify our thin $\mc{P}_{\nu^+}(\la)$-list with a thin $\mc{P}_{\nu^+}(\nu^+)$-list; an ineffable branch in one will correspond to an ineffable branch in the other.

To finish the proof, it suffices to verify the assumptions of Theorem \ref{thm:gen1cardITP} in $V[K]$. Since $\Coll(\nu^+, \la)$ is $\nu^+$-closed, each $\LL_\mu$ remains the product of a small and a closed forcing. Since $\ka_0$ is indestructibly supercompact in $V$, it will remain supercompact in $V[K]$.

By assumption, for each $\rho < \tau$, there is a generic $\la$-supercompactness embedding $j_{\rho+2}$ with domain $V$ and critical point $\ka_{\rho+2}$, added by a poset $\FF$ whose $\ka_\rho$-power is $<\ka_\rho$-distributive in $V$. We wish to extend the domain of $j_{\rho+2}$ to $V[K]$. Note that $\bigcup j''K$ is a condition in $j(\Coll(\nu^+,\la)^{V[K]}) = \Coll(j(\nu^+, j(\la))^{V[K]})$, so the embedding will be contained in $V[K][F][K^*]$, where $F$ is generic for $\FF$ and $K^*$ is a generic for $\Coll(j(\nu^+), j(\la))^{V[K]}$ containing $\bigcup j''K$. Noting that $\FF\times\Coll(j(\nu^+), j(\la))^{V[K]}$ will have a $<\ka_{\rho}$-distributive $\ka_\rho$-power, it follows that the second condition of Theorem \ref{thm:gen1cardITP} holds.
We conclude that there is $\mu \in \Index$ such that in $V[K][\LL_\mu]$, $d^\mu$ has an ineffable branch. This concludes the proof.

\begin{cor}\label{cor:anycof}
	Let $\tau$ be a regular cardinal with $\tau < \aleph_\tau$, and let $\langle \ka_\rho \mid \rho < \tau\rangle$ be a continuous increasing sequence of cardinals with supremum $\nu$ such that the following holds:
	\begin{itemize}
		\item $\ka_{\rho+1} = \ka_{\rho}^+$ for all limit ordinals $\rho < \tau$
		\item $\ka_n$ is indestructibly supercompact for all $n < \w$
		\item $\ka_{\rho+2}$ is indestructibly supercompact for all $\w\leq \rho < \tau$
	\end{itemize}
 	Then there is a generic extension in which ITP holds at $\aleph_{\tau+1}$, and $\tau$ remains a regular cardinal.
\end{cor}
\begin{proof}
	Define
	 \[\HH = \left(\prod _{\rho < \w} \Coll(\ka_\rho, <\ka_{\rho+1})\right)\times \left(\prod_{\w\leq \rho < \tau }\Coll(\ka_{\rho+1}, <\ka_{\rho+2})\right),\]
	 and let $H$ be generic for $\HH$. Let $\Index$ be the set of all $\mu < \ka_0$ of cofinality $\tau$. For all $\mu \in \Index$, let $\LL_{\mu} = \Coll(\tau, \mu) \times \Coll(\mu^{++}, <\ka_0)$.
	
	We apply Theorem \ref{thm:gen2cardITP} to $V[H]$ to obtain $\mu \in \Index$ such that in the extension of $V[H]$ by $\LL_\mu$, ITP holds at $\nu^+$. In this extension, $\ka_{\rho+1} = \ka_{\rho}^+$ for all $\rho < \tau$, while $\ka_0$ becomes $\tau^{+3}$, so $\nu = \aleph_\tau$ and $\nu^+ = \aleph_{\tau+1}$. Note also that since every poset is either $\tau$-cc or $\tau$-closed, $\tau$ remains a regular cardinal.
\end{proof}

\section{ITP from \texorpdfstring{$\aleph_4$}{ℵ\_4} to \texorpdfstring{$\aleph_{\w+1}$}{ℵ\_ω+1}}\label{s:ITPalmostuptoNw+1}

When we attempt to apply the results of the previous section to Neeman's construction, we run into several issues: in particular, $\LL_{\mu}$ will not be sufficiently closed, and Theorem \ref{thm:gen2cardITP} doesn't work when $i$ is a generic embedding. We can avoid these issues by making the following modifications:
\begin{itemize}
	\item We index our list of supercompact cardinals starting at $3$ rather than $2$. Our choice of $\mu \in \Index$ will select not only $\ka_1 \defeq \mu^+$ but also $\ka_2 \defeq \mu^{++}$.
	\item We replace $\A_0$ and $\A_1$ with the trivial forcings.
	\item We replace $\CC_0$ with the trivial poset, since we do not need any collapses to obtain $\ka_1^+ = \ka_2$.
	\item Since $\A_0$ and $\A_1$ are trivial, we can remove most of the restrictions on $\Index$.
	\item Since $\A_1$ is trivial, the poset $\CC_1$ will not collapse cardinals between $\ka_2$ and $\ka_3$, so we replace it with the trivial poset. We add $\Coll(\ka_2, <\ka_3)$ to the product. Note that this poset depends on the parameter $\ka_2 = \mu^{++}$.
	\item Since $\A_0, \A_1$, $\CC_0$, and $\CC_1$ are all trivial, we do not require the initial stages of Laver preparation, so we can set $\UU_0, \BB_0, \UU_1$, and $\BB_1$ to likewise be trivial.
	\item We replace $\A_2$ with $\sum_{\mu \in \Index} \Add(\mu^{++}, \ka_4)$.
\end{itemize}
Note that $\ka_2$ will no longer be generically supercompact, so the tree property (and thus also ITP) will not hold at $\aleph_2$ in the final model. Note also that since we removed $\A_1$, GCH will hold at $\aleph_1$, so the tree property will also fail at $\aleph_3$.

We now formally define the construction.
Let $\langle \ka_n \mid 3 \leq n < \w\rangle$ be an increasing sequence of indestructibly supercompact cardinals with supremum $\nu$, and suppose there is a partial function $\phi$ such that for all $n$, $\phi\uhr \ka_n$ is an indestructible Laver function for $\ka_n$ and for all $\alpha \in \dom(\phi)$, if $\gamma$ is in $\dom(\phi)\cap \alpha$ then $\phi(\gamma) \in V_\alpha$.

\begin{defn}
	We define the set Index as the set of all $\mu < \ka_3$ so that $\mu$ is a strong limit cardinal of cofinality $\w$ and $\dom(\phi)$ has a largest point $\la$ below $\mu$.
\end{defn}
Note that the elements of $\{\mu^+ \mid \mu \in \Index \}$ are the potential options for $\aleph_1$, and $\{\mu^{++} \mid \mu \in \Index\}$ are the potential options for $\aleph_2$. (Unlike in Section \ref{s:TPuptoNw+1}, we can define this set in advance because $A_0$ and $A_1$ are trivial, so Index won't depend on them.)

\begin{defn}
	For $n \geq 3$, let $\A_n \defeq \Add(\ka_n, \ka_{n+2})$. Let $\ka_0$ denote $\w$, and set $\A_0$ and $\A_1$ to be the trivial forcings. Let $A_2 = \sum_{\mu \in \text{Index}} \Add(\mu^{++}, \ka_4)$. We use $\ka_1$ and $\ka_2$ to refer to the values $\mu^+$ and $\mu^{++}$ chosen by a fixed generic. Let $\A$ be the full support product of $\A_n$ for $n < \w$.
\end{defn}

The Laver preparation posets $\BB$ and $\UU$ are constructed exactly as in Definition \ref{def:laver}
, except that we define $\BB_0$, $\BB_1$, $\UU_0$, and $\UU_1$ to be trivial. Since $\A_0$ and $\A_1$ are trivial, and we won't be using $\CC$ to collapse cardinals below $\ka_3$, we do not need any Laver preparation for those cardinals.
The collapsing poset $\CC_n$ is as defined in Definition \ref{def:collapses}, except that we define $\CC_0$ and $\CC_1$ to be trivial.

Let $A$ be generic for $\A$, $U$ be generic for $\UU$ over $V[A]$, $S$ be generic for $\CC^{+A\ast U}$ over $V[A][U]$, and $e$ be generic for $\Coll(\w, \mu)\times \Coll(\mu^{++}, <\ka_3)$ over $V[A][U][S]$, noting that $\mu$ is determined by the generic for $\A_2$ included in $A$.

The cardinal structure in $V[A][U][S][e]$ is similar to that in Neeman's construction.
\begin{claim}
	In $V[A][U][S][e]$, the following properties hold.
	\begin{itemize}
		\item $\ka_n = \aleph_n$ for each $n$.
		\item $2^{\ka_n} = \ka_{n+2}$ for each $n > 1$.
	\end{itemize}
\end{claim}
\begin{proof}
	Clearly $\ka_1 = \aleph_1$, and $\ka_2 = \ka_1^+$ in $V$. To verify that $\ka_2 = \aleph_2$, we need only check that $\ka_2$ is not collapsed. This follows from \cite[Claims 4.18 - 4.21]{NeemanTPNw+1}. Similarly, the cardinals between $\ka_2$ and $\ka_3$ are collapsed, while $\ka_3$ is preserved. For $n \neq 4$, the proof that $\ka_n = \aleph_n$ is identical to the proof of \cite[Lemma 4.24]{NeemanTPNw+1}. Since $\ka_1 = \aleph_1$, and $\ka_2 = \ka_1^+$ in $V$, we need only check that $\ka_2$ is not collapsed. This follows from \cite[Claims 4.18 - 4.21]{NeemanTPNw+1}.
	
	The proof that $2^{\ka_n} = \ka_{n+2}$ for $n > 1$ is identical to the proof of \cite[Claim 4.28]{NeemanTPNw+1}. (Note that this proof relies on several claims and lemmas - while our construction is slightly different, the proofs of those claims and lemmas are identical.)
\end{proof}

\begin{lemma}\label{lem:ITPbelow}
	In $V[A][U][S][e]$, ITP holds at $\aleph_{n+2}$ for all $n > 1$.
\end{lemma}
\begin{proof}
	The proof is exactly as in Subsection \ref{section:ITPuptoNw}. While the model has changed slightly, the proofs of the various facts used are identical to the corresponding proofs in \cite{NeemanTPNw+1} for $n \geq 2$.
\end{proof}

\begin{theorem}
	Let $\langle \ka_n \mid 3 \leq n < \w\rangle$ be an increasing sequence of supercompact cardinals with supremum $\nu$. Then there is a forcing extension in which $\aleph_n = \ka_n$, $\aleph_\w$ is strong limit, $(\nu^+)^V = \aleph_{\w+1}$, and ITP holds at $\aleph_{\w+1}$ and at $\aleph_n$ for $3 < n < \w$.
\end{theorem}
\begin{proof}
	As before, we define an intermediate model $M = V[A_{[3,\w)}][B_{[2,\w)}][C_{[2,\w)}]$, and let $\LL_\mu = \Add(\mu^{++}, \ka_4)^V \times \Coll(\mu^{++},<\ka_3) \times \Coll(\w,\mu)$. Note that $|\Coll(\w, \mu)| < \mu^+$, and $\Coll(\mu^{++}, <\ka_3)$ is $\mu^{++}$-closed.
	
	Next, we need to verify that the remainder of the forcing is $\mu^{++}$-closed in $M$. $\Add(\mu^{++}, \ka_4)^V$ is $\mu^{++}$-closed in $V$. To show that this poset is closed in $M$, we need to show that all of the intermediate posets are likewise closed. $\A_{[3,\w)}$ is $\ka_3$-closed. By \cite[Claim 4.7]{NeemanTPNw+1}, noting that $\A_0, \A_1, \UU_0,$ and $\UU_1$ are all trivial, we see that $\BB_{[2,\w)}$ is $\ka_2 = \mu^{++}$-closed over $V$. Finally, $\CC^{+B}\uhr[\ka_3, \nu)$ is $\ka_3$-closed over $V[B]$ by \cite[Claim 4.15]{NeemanTPNw+1}. We conclude that $\LL_\mu$ is the product of a $\mu^{++}$-closed poset with a poset of size $<\mu^+$.
	
	Since $\ka_3$ is supercompact, and the posets adding the embeddings for $n > 3$ have suitably distributive powers, we meet the hypotheses of Theorem \ref{thm:gen2cardITP}. It follows that there exists $\mu \in \Index$ such that in the extension of $M$ by $\LL_{\mu}$, ITP holds at $\aleph_{\w+1}$.
	
	Let $A_2 \times S_1 \times e$ be generic for $\LL_\mu$ over $M$, and let $A = A_2 \times A_{[3,\w)}$. Let $U = U_{[2,\w)}$ be the upwards closure of $B_{[2,\w)}$ in $\UU_{[2,\w)}$. Let $S_{[2,\w)}$ be the upwards closure of $C_{[2,\w)}$ in $\CC_{[1,\w)}^{+A\ast U}$, and let $S = S_1 \times S_{[2,\w)}$. Let $N$ be the model $V[A][U][S]$.
	
	\begin{claim}\label{cl:TPNeemanLem5.9}
		There is a $\mu^+$-closed extension $N[G]$ of $N$ such that $G$ is still generic over $N[e]$, both $\nu$ and $\nu^+$ are still cardinals in $N[e][G]$, and $M[A_2 \times S_1 \times e]$ is contained in $N[e][G]$.
	\end{claim}
	\begin{proof}
		Analogous to Proposition \ref{prop:TPNeemanLem5.9}. The relevant forcing is the product of the factor forcing refining $U_{[2,\w)}$ to a filter for $\BB^{+A_2\ast U_2}\uhr[\ka_3, \nu)$ and the factor forcing refining $S_{[2,\w)}$ to a generic for $\CC^{+A_2\ast U_2}\uhr[\ka_3, \nu)$. These posets are in $V[A][U][S]$. That they are $\mu^+$-closed in $V[A][U][S]$ follows from the proofs of \cite[Claim 4.9]{NeemanTPNw+1}, \cite[Claim 4.16]{NeemanTPNw+1}, and \cite[Claim 4.26]{NeemanTPNw+1}; preservation of cardinals follows from \cite[Claims 4.18-4.21]{NeemanTPNw+1}.
	\end{proof}
	
	\begin{claim}
		In $N[e]$, $ITP$ holds at $\nu^+$.
	\end{claim}
	\begin{proof}
		Identical to the proof of Claim \ref{cl:TPnu^+inN[e]}, using Claim \ref{cl:TPNeemanLem5.9} instead of Proposition \ref{prop:TPNeemanLem5.9}.
	\end{proof}
	
	\begin{claim}
		In $N[e]$, $ITP$ holds at $\aleph_n$ for all $n \geq 4$.
	\end{claim}
	\begin{proof}
		Immediate from Lemma \ref{lem:ITPbelow}.
	\end{proof}
	
	We conclude that in $N[e]$, $\ITP$ holds at $\aleph_{\w+1}$ and at $\aleph_{n}$ for all $3 < n < \w$, completing the proof.
\end{proof}

\section{Successors of singular cardinals of multiple cofinalities}\label{s:manycofs}
Given any regular cardinal $\tau$, along with $\tau$-many supercompact cardinals, we can apply the techniques of the previous sections to obtain the strong or super tree properties at $\aleph_{{\tau}+1}$. If we try to apply these techniques for multiple cofinalities $\tau_0$ and $\tau_1$ simultaneously, there is a major obstacle: in the resulting model, $\tau_1$ will be collapsed to $\tau_0$, so $\aleph_{\tau_1}$ will only have cofinality $\tau_0$. Instead, if $\tau_1$ is of the form $\omega_{\alpha_1}$ for some successor ordinal $\alpha_1$, we can ensure that $\aleph_{\omega_{\alpha_1}+1}$ \emph{as computed in the final model} has ITP. In fact, this can be done for finitely many such cardinals.

In this section, given any finite sequence $\alpha_0, \dots, \alpha_n$ of ordinals where $\alpha_i < \omega_{\alpha_i}$ and $\aleph_{\alpha_i}$ is a regular cardinal, we construct a model where the strong tree property or ITP hold at each $\aleph_{\omega_{\alpha_n}+1}$ simultaneously. First we demonstrate the construction for $\aleph_{\w+1}$ and $\aleph_{\w_1+1}$, and then we present the construction in full generality. Note that the large cardinal hypotheses are somewhat stronger than might be expected; we require supercompact-many supercompacts, instead of $\w_1$-many supercompacts. This is because we want to obtain ITP at $\aleph_{\w_1+1}$ in our final model $V[H][\LL_{\mu,\delta}]$, which will be defined below. In this model, $\w_1 = \mu^+$, so we require $\mu^+$-many supercompact cardinals; since $\mu$ is not fixed ahead of time, and the only upper bound we have is $\ka_0$, we need $\ka_0$-many supercompacts. The same argument, using only $\w_1$-many supercompacts, can obtain ITP at $\aleph_{(\w_1)^V+1}$, but since $(\w_1)^V$ is collapsed by $\LL_\mu$ (for all $\mu > \w$) this will be the successor of a singular of only countable cofinality.

\begin{thm}\label{thm:tpw+1andw_1+1}
	Let $\ka_0$ be indestructibly supercompact, and let $\langle \ka_\rho \mid \rho < \ka_0 \rangle$ be an increasing continuous sequence of cardinals such that the following holds:
	\begin{itemize}
		\item $\ka_{\rho+1} = \ka_{\rho}^+$ for all limit ordinals $\rho < \ka_0$
		\item $\ka_n$ is indestructibly supercompact for all $n < \w$
		\item $\ka_{\rho+2}$ is indestructibly supercompact for all $\w\leq \rho < \ka_0$.
	\end{itemize}
	 Then there is a generic extension in which ITP holds at $\aleph_{\w_1+1}$ and the strong tree property holds at $\aleph_{\w+1}$.	
\end{thm}
\begin{proof}
	As in Corollary \ref{cor:anycof}, we define
	\[\HH = \left(\prod_{\rho<\w} \Coll(\ka_\rho, <\ka_{\rho+1})\right) \times \left(\prod_{\w\leq\rho<\ka_0} \Coll(\ka_{\rho+1}, <\ka_{\rho+2})\right).\]
	
	Let $H$ be generic for $\HH$.
	Let $\Index = \{(\mu, \delta) \mid \cf(\mu) = \w, \cf(\delta) = \mu^+, \mu < \delta < \ka_0\}.$ For each pair $(\mu, \delta) \in \Index$, let $\LL_{\mu,\delta} = \Coll(\w, \mu) \times \Coll(\mu^+, \delta) \times \Coll(\delta^{++}, <\ka_0)$. We wish to find $(\mu, \delta) \in \Index$ such that in the generic extension of $V[H]$ by $\LL_{\mu,\delta}$, ITP holds at $\aleph_{\w_1+1}$ and the strong tree property holds at $\aleph_{\w+1}$.
	
	Let $I$ be the projection of $\Index$ to the first coordinate. For all $\mu \in I$, let $I_\mu = \{\delta \mid \mu < \delta < \ka_0, \cf(\delta) = \mu^+\}$.
	For each $\mu \in I$, consider the initial segment $\langle \ka_\rho \mid \rho < \mu^+ \rangle$, with supremum $\nu_\mu \defeq \ka_{\mu^+}$. We now verify the hypotheses of Theorem \ref{thm:gen2cardITP}.
	 
	For any $\delta \in I_\mu$, $\Coll(\w, \mu) \times \Coll(\mu^+, \delta)$ has size $< \delta^+$, while $\Coll(\delta^{++}, <\ka_0)$ is $\delta^{++}$-closed; moreover $|\LL_{\mu,\delta}|<\ka_2$.
	
	In $V[H]$, $\ka_0$ is $\nu_\mu^+$-supercompact; let $i$ be the corresponding embedding. Noting that $[x \mapsto \ka_x^{+\mu^+}]_{U_0} = \nu_\mu$, and $\ka_x^{+\mu^+} \in I_\mu$, we conclude that $\nu_\mu\in i(I_\mu)$. Note also that there are generic supercompactness embeddings $j_{\rho+2}$ added by $\Coll(\ka_{\rho+1}, j(\ka_{\rho+3}))$. This poset is $\ka_{\rho+1}$-closed in $V[H\uhr[\rho+2, \ka_0)]$, and $V[H]$ is a $\ka_{\rho+1}$-cc extension of this poset, so in particular the product $\Coll(\ka_{\rho+1}, j(\ka_{\rho+3}))^{\ka_\rho}$ is $<\ka_\rho$-distributive.
	
	Thus for each $\mu \in I$, we can apply Theorem \ref{thm:gen2cardITP} with respect to this initial segment and $I_\mu$ to conclude that there exists some $\delta_\mu \in I_\mu$ such that ITP holds at $\nu_\mu^+$ in the generic extension of $V[H]$ by $\LL_{\mu,\delta_\mu}$.
	
	We will now pin down the first coordinate. Consider the initial segment $\langle \ka_i \mid i<\w\rangle$, with supremum $\nu \defeq \ka_\w$.
	In $V[H]$, $\ka_0$ is $\nu^+$-supercompact; let $i$ be the corresponding embedding. Noting that $[x \mapsto \ka_x^{+\w}]_{U_0} = \nu$, and $\ka_x^{+\w}$ is in $I$, we conclude that $\nu \in i(I)$. Note also that $|\LL_{\mu, \delta_\mu}| < \ka_2$ for all $\mu \in I$, so we also meet the hypothesis of Theorem \ref{thm:TP2cardfullNA} (with respect to the first $\w$-many supercompacts, using $I$ as our index set). Thus there exists $\mu \in I$ such that the strong tree property holds at $\nu^+$ in the extension of $V[H]$ by $\LL_{\mu, \delta_\mu}$. Note that the pair $(\mu, \delta_\mu)$ is in $\Index$. Then in the extension of $V[H]$ by a generic $L$ for $\LL_{\mu, \delta_\mu}$, ITP will hold at $\nu_\mu^+$, and the strong tree property will hold at $\nu^+$.
	
	To finish the proof, we examine the cardinal structure in the new model. In $V[H][L]$, $\mu$ is collapsed to $\w$, so $\mu^+ = \aleph_1$. The other parameter $\delta$ is collapsed to $\mu^+$, so $\delta^+ = \mu^{++} = \aleph_2$. After forcing with $\HH$, we have that $\ka_{\rho^+} = \ka_{\rho+1}$ for all $\rho < \nu$, and in $V[H][L]$, $\ka_0$ becomes $\aleph_4$. We conclude that since $\ka_\w^+$ is the successor of the limit of $\langle \ka_n \mid n < \w\rangle$, $\nu^+ = \aleph_{w+1}$. Similarly, $\nu_\mu^+$ is the successor of the limit of $\langle \ka_\rho \mid \rho < \w_1\rangle$, so $\nu_\mu^+ = \aleph_{\w_1+1}$.
\end{proof}
Now we generalize this argument to apply to any finite sequence of cofinalities.

\begin{thm}
	Let $\alpha_0, \dots, \alpha_n$ be an increasing sequence of ordinals, such that for all $i\leq n$, $\omega_{\alpha_i}$ is a regular cardinal with $\alpha_i < \omega_{\alpha_i}$.
	Let $\ka_0$ be indestructibly supercompact, and let $\langle \ka_\rho \mid \rho < \ka_0 \rangle$ be an increasing continuous sequence of cardinals such that the following holds:
	\begin{itemize}
		\item $\ka_{\rho+1} = \ka_{\rho}^+$ for all limit ordinals $\rho < \ka_0$
		\item $\ka_n$ is indestructibly supercompact for all $n < \w$
		\item $\ka_{\rho+2}$ is indestructibly supercompact for all $\w\leq \rho < \ka_0.$
	\end{itemize}
	Then there is a generic extension in which the strong tree property holds at $\aleph_{\omega_{\alpha_i}+1}$ for all $i \leq n$. In addition, if $\alpha_{i+1} > \alpha_i+1$, then ITP holds at $\aleph_{\omega_{\alpha_i}+1}$.
\end{thm}
\begin{proof}
We define
\[\HH = \left(\prod_{\rho<\w} \Coll(\ka_\rho, <\ka_{\rho+1})\right) \times \left(\prod_{\w\leq\rho<\ka_0} \Coll(\ka_{\rho+1}, <\ka_{\rho+2})\right).\]
Let $H$ be generic for $\HH$. For all $0 \leq i < n$, define $f(i)$ as the unique ordinal $\beta$ such that $\alpha_{i+1} = \alpha_{i}+\beta$. In particular, this means that $\omega_{\alpha_{i+1}} = \omega_{\alpha_i}^{+f(i)}$. Let $\Index$ be the set of all increasing sequences $s = \langle \delta_i \mid i \leq n\rangle$ such that $\delta_n < \ka_0$, $\cf(\delta_0) = \omega_{\alpha_0}$, and for all $0 \leq i < n$, $\cf(\delta_{i+1}) = \delta_i^{+f(i)}$. For all $s \in \Index$, in $V$ we define
\[\LL_s \defeq \Coll(\omega_{\alpha_0}, \delta_0) \times \left(\prod_{1\leq i < n} \Coll(\cf(\delta_{i}), \delta_{i})\right) \times \Coll(\delta_n^{++}, <\ka_0).\]
Note that for all $s \in \Index$, $|\LL_s| < \ka_2$.

We wish to show that there exists $s \in \Index$ so that in the extension of $V[H]$ by $\LL_s$, the strong tree property holds at $\aleph_{\omega_{\alpha_i}+1}$ for all $i \leq n$, and if $\alpha_{i+1} > (\alpha_i)+1$, then ITP holds at $\aleph_{\omega_{\alpha_i}+1}$.
	
First, we note that in $V[H]$, for any $\la \geq \nu^+$, there are generic $\la$-supercompactness embeddings $j_{\rho+2}$ with critical point $\ka_{\rho+2}$ added by $\Coll(\ka_{\rho+1}, j(\ka_{\rho+3}))$. This poset is $\ka_{\rho+1}$-closed in $V[H\uhr[\rho+2, \ka_0)]$, and $V[H]$ is a $\ka_{\rho+1}$-cc extension of this poset, so in particular the full support product $\Coll(\ka_{\rho+1}, j(\ka_{\rho+3}))^{\ka_\rho}$ is $<\ka_\rho$-distributive in $V[H]$. Note also that in $V[H]$, $\ka_0$ is supercompact. Let $\pi$ be a $\nu^+$-supercompactness embedding with critical point $\ka$.

We need to choose $\delta_i$ for all $i$; we do so inductively, beginning at $\delta_n$ and working downwards.

For every fixed $s = \langle \delta_0, \dots, \delta_{n-1}\rangle$ in the projection of $\Index$ to the first $n$ components, let $\delta_n^* = \delta_{n-1}^{+f(n-1)}$, and let $I_s = \{\delta < \ka_0 \mid \cf(\delta) = \delta^*_n\}$. Note that for all $\delta \in I_{s}$, $s\cat \delta \in \Index$. Consider the initial segment $\langle \ka_\rho \mid \rho < \delta_n^*\rangle$, and let $\nu_s^n$ be the supremum of this sequence.

Note that for each $\delta \in I_{s}$, $\Coll(\delta_n^{++}, <\ka_0)$ is $\delta^{++}$-closed, and the remainder of $\LL_{s\cat \delta}$ has size $< \delta^+$. Note also that $\ka_x^{+\delta^*}$ is in $I_{s}$ for all $x$, so $\nu_s^n \in \pi(I_s)$. Then applying Theorem \ref{thm:gen2cardITP}, we conclude that for each sequence $s$, there exists some $\delta^s_n \in I_{s}$ such that in the extension of $V[H]$ by $\LL_{s\cat\delta^s_n}$, ITP holds at $(\nu_s^n)^+$.

Now let $0 < k  < n$, and let $s$ be a fixed sequence $\langle \delta_0, \dots, \delta_{k-1}\rangle$ in the projection of $\Index$ to the first $k$ components. Let $\nu_s^i = \sup \langle \ka_\rho \mid \rho < \cf(\delta_i)\rangle$. Let $\delta_k^* = \delta_{k-1}^{+f(k-1)}$ and let $I_s = \{\delta < \ka_0 \mid \cf(\delta) = \delta_k^*\}.$ Working inductively, assume that for all $\delta \in I_s$ we have determined a sequence $s'_\delta = \langle \delta_i \mid k < i \leq n\rangle$ of length $n-k$ such that the sequence $s_\delta = s \cat\delta\cat s'_\delta$ is in $\Index$, and in the extension of $V[H]$ by $\LL_{s_\delta}$, ITP or the strong tree property holds at $\nu_{s_\delta}^i$ for all $k < i \leq n$.

Consider the initial segment $\langle \ka_\rho \mid \rho < \delta_{k-1}^{+f(k)}\rangle$, and let $\nu_s^k$ be the supremum of this sequence. Note that $\ka_x^{+\delta_k^*}$ is in $I_{s}$ for all $x$, so $\nu_s^n \in i(I_s)$.

We factor $\LL_{s_\delta}$ as $\LL_s^0 \times \LL_s^1$, where each term is defined as follows, with $\delta_k$ denoting $\delta$:
\begin{align*}
\LL_{s_\delta}^0 &\defeq \Coll(\w_{\alpha_0}, \delta_0) \times \prod_{1 \leq i \leq k} \Coll(\cf(\delta_{i}), \delta_{i}), \\
\LL_{s_\delta}^1 &\defeq \left(\prod_{k < i < n } \Coll(\cf(\delta_{i}), \delta_{i})\right) \times \Coll(\delta_n^{++}, <\ka_0).
\end{align*}

Note that for all $\delta \in I_s$, $|\LL_{s_\delta}^0| < \delta^+$. If $\alpha_i > \alpha_{i}+1$, then $\cf(\delta_{k+1}) > \delta^+$, so we see that $\LL_{s_\delta}^1$ is $\delta^{++}$-closed. We apply Theorem \ref{thm:gen2cardITP} to choose some $\delta_k \in I_s$ such that in the extension of $V[H]$ by $\LL_{s,\delta_k}$, ITP holds at $(\nu_s^k)^+$.

If $\cf(\delta_{k+1}) = \delta^+$, we will not have enough closure to  apply Theorem \ref{thm:gen2cardITP}. However, noting that $|\LL_{s_\delta}|< \ka_1$, we can instead apply Theorem \ref{thm:TP2cardfullNA} to select $\delta_k$ such that in the extension of $V[H]$ by $\LL_{s_{\delta_k}}$, the strong tree property holds at $(\nu_s^k)^+$.

We repeat this argument for $k=0$. Let $\nu$ be the limit of $\langle\ka_\rho \mid \rho <\omega_{\alpha_0}\rangle$. By the previous stages in our construction, for each $\mu$ in the projection of $\Index$ to the first coordinate, we have built a sequence $s_\mu = \mu \cat \langle \delta_i^\mu \mid 0 < i \leq n\rangle$ in $\Index$ such that in the extension of $V[H]$ by $\LL_{s_\mu}$, $\ITP$ or the strong tree property holds at $\nu_{s_\mu}^i$ for all $i > 0$. If $\cf(\delta_1^\mu) > \mu^+$, then setting $\LL_{s_\mu}^0 = \Coll(\omega_{\alpha_0}, \delta_0)$ and $\LL_{s_\mu}^1$ to be the remainder of the forcing, we can apply Theorem \ref{thm:gen2cardITP} to obtain some $\mu$ such that after forcing with $\LL_{s_\mu},$ ITP holds at $\nu^+$. If this is not the case, we can apply Theorem \ref{thm:TP2cardfullNA} to obtain the strong tree property at $\nu^+$. Let $\delta_0 = \mu$, and let $s \defeq s_\mu$.
By construction, $s \in \Index$ such that in the extension of $V[H]$ by $\LL_s$, ITP or the strong tree property hold at $\nu_i^+$ for all $i \leq n$. 

Finally, we examine the cardinal structure in this extension. The first chosen cardinal $\delta_0$ is collapsed to $\omega_{\alpha_0}$, while $\delta_1$ is collapsed to $\delta_0^{+f(1)} = \omega_{\alpha_0}^{+f(1)} = \omega_{\alpha_1}$. Working inductively we see that each $\delta_i$ is collapsed to $\omega_{\alpha_{i}}$. The cardinals between each successive $\ka_\rho$ are collapsed, and $\ka_0$ will become $\omega_{\alpha_n}^{+3}$. It follows that each $\nu_i$ becomes $\ka_0^{+\omega_{\alpha_i}} = \aleph_{\omega_{\alpha_i}}$.

We conclude that, in the extension of $V[H]$ by $\LL_s$, the strong tree property holds at $\aleph_{\omega_{\alpha_i}+1}$ for all $i \leq n$, and if $\alpha_{i+1} > \alpha_i+1$, ITP holds at $\aleph_{\omega_{\alpha_i}+1}$.
\end{proof}

\section{Open Problems}\label{s:open}

The natural question is, of course, can we get ITP at $\aleph_{\w+1}$ along with ITP at every $\aleph_n$? As we have shown, the difficulty is at $n = 2$ and $n=3$.

\begin{quest}
	Can we obtain ITP at $\aleph_{\w+1}$ and at $\aleph_2$ simultaneously? What about at $\aleph_{\w+1}$ along with $\aleph_n$ for all $2 \leq n < \w?$
\end{quest}

Another more ambitious project is obtaining the tree property, or more optimistically its generalizations, up to $\aleph_{\w_1+1}$. There are two major obstacles. First, we would need to obtain the tree property at the successor of every cardinal with countable cofinality. Golshani and Hayut have shown in \cite{golshani-hayut:tpcountablesegment} that this can be done for any countable initial segment of these cardinals. (Using the lemmas in this paper, their argument immediately generalizes to the one-cardinal version of ITP, but not to the full two-cardinal property.)
\begin{quest}
Can we obtain the tree property (or ITP) at every successor of a singular cardinal below $\aleph_{\w_1}$?
\end{quest}

The other obstacle comes from the reflection principles used to obtain the tree property at the successor of a singular cardinal. To obtain the tree property at $\aleph_{\w_1+1}$, $\aleph_2$ in the final model will be $\mu^+$, where $\mu$ is the chosen cardinal that is collapsed to $\w_1$. To obtain the tree property at $\aleph_2$ using standard techniques, however, $\aleph_2$ must be generically supercompact.

\begin{quest}
	Can we obtain the tree property at $\aleph_2$ and at $\aleph_{\w_1+1}$ simultaneously? What about the strong tree property or ITP?
\end{quest}

Finally, while we have shown that it is consistent for ITP to hold simultaneously at the successors of small singular cardinals with different cofinalities, our argument is limited to finitely many cofinalities.

\begin{quest}
	Given a sequence $\langle \gamma_n \mid n < \w\rangle$ of regular cardinals, can we obtain the strong tree property (or ITP) at $\aleph_{\gamma_n+1}$ for all $n$ simultaneously?
\end{quest}

\bibliography{bib}
\bibliographystyle{plain}

\end{document}